%%%%%%%%%%%%%%%%%%  tex macros for preprints, cm version %%%%%%%%%%%%%%
%         (P. Ginsparg <ginsparg@lanl.gov>, last updated 7/94)
%         hypertex extensions (still provisional), 7/26/94
%	  Some modifications by C.R.Mafra, 2012

%comment out this line to restore non-hyper functionality
%\input hyperbasics

%\input amssym.tex % for blackboard bold
\input epsf.tex

\def\unredoffs{}
\tolerance=1000\hfuzz=2pt
\catcode`\@=11 % This allows us to modify PLAIN macros.
\ifx\hyperdef\UNd@FiNeD\def\hyperdef#1#2#3#4{#4}\def\hyperref#1#2#3#4{#4}\def\href#1#2{#2}\fi
\magnification=1200\unredoffs\baselineskip=16pt plus 2pt minus 1pt
%
%
% (restores pagenumbers)

%

%%%%%% Hour:Minute %%%%%%%%%%%%%%%%%
{\count255=\time\divide\count255 by 60 \xdef\hourmin{\number\count255}
 \multiply\count255 by-60\advance\count255 by\time
 \xdef\hourmin{\hourmin:\ifnum\count255<10 0\fi\the\count255}
}
\def\date{\number\day.\number\month.\number\year\ at \hourmin}

%%%%%%%%%%%% Draft mode %%%%%%%%%%%%%
% puts date/time on each page in big mode, writes labels in margins

% use \nolabels to get rid of eqn, ref, and fig labels in draft mode
\def\nolabels{\def\wrlabeL##1{}\def\eqlabeL##1{}\def\reflabeL##1{}}
\def\writelabels{\def\wrlabeL##1{\leavevmode\vadjust{\rlap{\smash%
{\line{{\escapechar=` \hfill\rlap{\sevenrm\hskip.03in\string##1}}}}}}}%
\def\eqlabeL##1{{\escapechar-1\rlap{\sevenrm\hskip.05in\string##1}}}%
\def\reflabeL##1{\noexpand\llap{\noexpand\sevenrm\string\string\string##1}}}
\nolabels

% tagged sec numbers
\global\newcount\secno \global\secno=0
\global\newcount\meqno \global\meqno=1
\def\s@csym{}

%%%%%%%%% Section %%%%%%%%%%%%%
\def\newsec#1\par{\global\advance\secno by1%
{\toks0{#1}\message{(\the\secno. \the\toks0)}}%
\global\subsecno=0\eqnres@t\let\s@csym\secsym\xdef\secn@m{\the\secno}\noindent
{\bf\hyperdef\hypernoname{section}{\the\secno}{\the\secno.} #1}%
\writetoca{{\string\hyperref{}{section}{\the\secno}{\bf \the\secno\quad}} {\bf #1}}\par%
\nobreak\medskip\nobreak\noindent\ignorespaces}
\def\eqnres@t{\xdef\secsym{\the\secno.}\global\meqno=1\bigbreak\bigskip}
\def\sequentialequations{\def\eqnres@t{\bigbreak}}\xdef\secsym{}

%%%%%%%% Subsection %%%%%%%%%%%
\global\newcount\subsecno \global\subsecno=0
\def\subsec#1\par{\global\advance\subsecno by1%
{\toks0{#1}\message{(\s@csym\the\subsecno. \the\toks0)}}%
\global\subsubsecno=0%
\ifnum\lastpenalty>9000\else\bigbreak\fi
\noindent{\it\hyperdef\hypernoname{subsection}{\secn@m.\the\subsecno}%
{\secn@m.\the\subsecno.} #1}\writetoca{\string\hskip1.45cm
{\string\hyperref{}{subsection}{\secn@m.\the\subsecno}{\secn@m.\the\subsecno.}}
{#1}}\par\nobreak\medskip\nobreak\noindent\ignorespaces}

%%%%%%% Appendix %%%%%%%%%%%%%%
\def\appendix#1#2{\global\meqno=1\global\subsecno=0\xdef\secsym{\hbox{#1.}}%
\bigbreak\bigskip\noindent{\bf Appendix \hyperdef\hypernoname{appendix}{#1}%
{#1.} #2}{\toks0{(#1. #2)}\message{\the\toks0}}%
\xdef\s@csym{#1.}\xdef\secn@m{#1}%
\writetoca{{\string\hyperref{}{appendix}{#1}{\bf {#1}\quad}} {\bf #2}}%
\par\nobreak\medskip\nobreak}

% \eqn\label{a+b=c}   gives displayed equation, numbered consecutively within sections.
% \eqnn, \eqna        define labels in advance, use \eqna\label before an eqalign and
%                     later \label a, \label b etc inside eqalign to get (2.3a), (2.3b) etc
%
\def\checkm@de#1#2{\ifmmode{\def\f@rst##1{##1}\hyperdef\hypernoname{equation}%
{#1}{#2}}\else\hyperref{}{equation}{#1}{#2}\fi}
\def\eqnn#1{\DefWarn#1\xdef #1{(\noexpand\relax\noexpand\checkm@de%
{\s@csym\the\meqno}{\secsym\the\meqno})}%
\wrlabeL#1\writedef{#1\leftbracket#1}\global\advance\meqno by1}
\def\f@rst#1{\c@t#1a\em@ark}\def\c@t#1#2\em@ark{#1}
\def\eqna#1{\DefWarn#1\wrlabeL{#1$\{\}$}%
\xdef #1##1{(\noexpand\relax\noexpand\checkm@de%
{\s@csym\the\meqno\noexpand\f@rst{##1}1}{\hbox{$\secsym\the\meqno##1$}})}
\writedef{#1\numbersign1\leftbracket#1{\numbersign1}}\global\advance\meqno by1}
\def\eqn#1#2{\DefWarn#1%
\xdef #1{(\noexpand\hyperref{}{equation}{\s@csym\the\meqno}%
{\secsym\the\meqno})}$$#2\eqno(\hyperdef\hypernoname{equation}%
{\s@csym\the\meqno}{\secsym\the\meqno})\eqlabeL#1$$%
\writedef{#1\leftbracket#1}\global\advance\meqno by1}
\def\xeqn{\expandafter\xe@n}\def\xe@n(#1){#1}
\def\xeqna#1{\expandafter\xe@n#1}
\def\eqns#1{(\e@ns #1{\hbox{}})}
\def\e@ns#1{\ifx\UNd@FiNeD#1\message{eqnlabel \string#1 is undefined.}%
\xdef#1{(?.?)}\fi{\let\hyperref=\relax\xdef\next{#1}}%
\ifx\next\em@rk\def\next{}\else%
\ifx\next#1\xeqn#1\else\def\n@xt{#1}\ifx\n@xt\next#1\else\xeqna#1\fi
\fi\let\next=\e@ns\fi\next}

\def\DefWarn#1{\ifx\UNd@FiNeD#1\else
\immediate\write16{*** WARNING: the label \string#1 is already defined ***}\fi}
%
% footnotes
\newskip\footskip\footskip14pt plus 1pt minus 1pt %sets footnote baselineskip
\def\footnotefont{\ninepoint}\def\f@t#1{\footnotefont #1\@foot}
\def\f@@t{\baselineskip\footskip\bgroup\footnotefont\aftergroup\@foot\let\next}
\setbox\strutbox=\hbox{\vrule height9.5pt depth4.5pt width0pt}
\global\newcount\ftno \global\ftno=0
\def\foot{\global\advance\ftno by1\def\foot@rg{\hyperref{}{footnote}%
{\the\ftno}{\the\ftno}\xdef\foot@rg{\noexpand\hyperdef\noexpand\hypernoname%
{footnote}{\the\ftno}{\the\ftno}}}\footnote{$^{\foot@rg}$}}
%
%
%     \ref\label{text}
% generates a number, assigns it to \label, generates an entry.
% To list the refs on a separate page,  \listrefs
%
\global\newcount\refno \global\refno=1
\newwrite\rfile
\def\ref{[\hyperref{}{reference}{\the\refno}{\the\refno}]\nref}
\def\nref#1{\DefWarn#1%
\xdef#1{[\noexpand\hyperref{}{reference}{\the\refno}{\the\refno}]}%
\writedef{#1\leftbracket#1}%
\ifnum\refno=1\immediate\openout\rfile=\jobname.refs\fi
\chardef\wfile=\rfile\immediate\write\rfile{\noexpand\item{[\noexpand\hyperdef%
\noexpand\hypernoname{reference}{\the\refno}{\the\refno}]\ }%
\reflabeL{#1\hskip.31in}\pctsign}\global\advance\refno by1\findarg}
%	horrible hack to sidestep tex \write limitation
\def\findarg#1#{\begingroup\obeylines\newlinechar=`\^^M\pass@rg}
{\obeylines\gdef\pass@rg#1{\writ@line\relax #1^^M\hbox{}^^M}%
\gdef\writ@line#1^^M{\expandafter\toks0\expandafter{\striprel@x #1}%
\edef\next{\the\toks0}\ifx\next\em@rk\let\next=\endgroup\else\ifx\next\empty%
\else\immediate\write\wfile{\the\toks0}\fi\let\next=\writ@line\fi\next\relax}}
\def\striprel@x#1{} \def\em@rk{\hbox{}}
\def\lref{\begingroup\obeylines\lr@f}
\def\lr@f#1#2{\DefWarn#1\gdef#1{\let#1=\UNd@FiNeD\ref#1{#2}}\endgroup\unskip}
\def\semi{;\hfil\break}
\def\addref#1{\immediate\write\rfile{\noexpand\item{}#1}} %now unnecessary
\def\listrefs{\vfill\supereject\immediate\closeout\rfile\writestoppt
\baselineskip=\footskip\centerline{{\bf References}}\bigskip{\parindent=20pt%
\frenchspacing\escapechar=` \input \jobname.refs\vfill\eject}\nonfrenchspacing}
\def\startrefs#1{\immediate\openout\rfile=\jobname.refs\refno=#1}
\def\xref{\expandafter\xr@f}\def\xr@f[#1]{#1}
\def\refs#1{\count255=1[\r@fs #1{\hbox{}}]}
\def\r@fs#1{\ifx\UNd@FiNeD#1\message{reflabel \string#1 is undefined.}%
\nref#1{need to supply reference \string#1.}\fi%
\vphantom{\hphantom{#1}}{\let\hyperref=\relax\xdef\next{#1}}%
\ifx\next\em@rk\def\next{}%
\else\ifx\next#1\ifodd\count255\relax\xref#1\count255=0\fi%
\else#1\count255=1\fi\let\next=\r@fs\fi\next}
%

%
% this is ugly, but moore insists
\newwrite\ffile\global\newcount\figno \global\figno=1
\def\fig{fig.~\hyperref{}{figure}{\the\figno}{\the\figno}\nfig}
\def\nfig#1{\DefWarn#1%
\xdef#1{fig.~\noexpand\hyperref{}{figure}{\the\figno}{\the\figno}}%
\writedef{#1\leftbracket fig.\noexpand~\xfig#1}%
\ifnum\figno=1\immediate\openout\ffile=\jobname.figs\fi\chardef\wfile=\ffile%
{\let\hyperref=\relax
\immediate\write\ffile{\noexpand\medskip\noexpand\item{Fig.\ %
\noexpand\hyperdef\noexpand\hypernoname{figure}{\the\figno}{\the\figno}. }
\reflabeL{#1\hskip.55in}\pctsign}}\global\advance\figno by1\findarg}
\def\xfig{\expandafter\xf@g}\def\xf@g fig.\penalty\@M\ {}
\def\figs#1{figs.~\f@gs #1{\hbox{}}}
\def\f@gs#1{{\let\hyperref=\relax\xdef\next{#1}}\ifx\next\em@rk\def\next{}\else
\ifx\next#1\xfig #1\else#1\fi\let\next=\f@gs\fi\next}
%
%% because TeXlive 2011 is buggy wrt to tikz pictures with plain TeX..
\def\figin{\epsfcheck\figin}\def\figins{\epsfcheck\figins}
\def\epsfcheck{\ifx\epsfbox\UnDeFiNeD
\message{(NO epsf.tex, FIGURES WILL BE IGNORED)}
\gdef\figin##1{\vskip2in}\gdef\figins##1{\hskip.5in}% blank space instead
\else\message{(FIGURES WILL BE INCLUDED)}%
\gdef\figin##1{##1}\gdef\figins##1{##1}\fi}
\def\DefWarn#1{}
\def\figinsert{\goodbreak\topinsert}
\def\ifig#1#2#3{\DefWarn#1\xdef#1{fig.~\the\figno}
\writedef{#1\leftbracket fig.\noexpand~\the\figno}%
\figinsert\figin{\centerline{#3}}
\smallskip
\leftskip=0pt \rightskip=0pt
\baselineskip12pt\noindent
{{\bf Fig.~\the\figno}\ \ninepoint #2}
\medskip
\global\advance\figno by1\par\endinsert}
%%%%%%%%%%%%%%%%%%%%%%%%%%%%%%%%%%%%%%%%%%%%%%%%%%%%%%%%%
\newwrite\lfile
{\escapechar-1\xdef\pctsign{\string\%}\xdef\leftbracket{\string\{}
\xdef\rightbracket{\string\}}\xdef\numbersign{\string\#}}
\def\writedefs{\immediate\openout\lfile=label.defs \def\writedef##1{%
{\let\hyperref=\relax\let\hyperdef=\relax\let\hypernoname=\relax
 \immediate\write\lfile{\string\def\string##1\rightbracket}}}}%
\def\writestop{\def\writestoppt{\immediate\write\lfile{\string\pageno
 \the\pageno\string\startrefs\leftbracket\the\refno\rightbracket
 \string\def\string\secsym\leftbracket\secsym\rightbracket
 \string\secno\the\secno\string\meqno\the\meqno}\immediate\closeout\lfile}}
\def\writestoppt{}\def\writedef#1{}

% Section, subsection and appendix labels %
% Note that there must be a blanck line after \newsec,\subsec and before \seclab,\subseclab!
\def\seclab#1\par{\DefWarn#1%
\xdef #1{\noexpand\hyperref{}{section}{\the\secno}{\the\secno}}%
\writedef{#1\leftbracket#1}\wrlabeL{#1=#1}\par%
\nobreak\medskip\nobreak\noindent\ignorespaces}
\def\subseclab#1\par{\DefWarn#1%
\xdef #1{\noexpand\hyperref{}{subsection}{\the\secno.\the\subsecno}%
{\the\secno.\the\subsecno}}\writedef{#1\leftbracket#1}\wrlabeL{#1=#1}\par%
\nobreak\medskip\nobreak\noindent\ignorespaces}
\def\subsubseclab#1\par{\DefWarn#1%
\xdef#1{\noexpand\hyperref{}{subsubsection}{\the\secno.\the\subsecno.\the\subsubsecno}%
{\the\secno.\the\subsecno.\the\subsubsecno}}\writedef{#1\leftbracket#1}\wrlabeL{#1=#1}\par%
\nobreak\medskip\nobreak\noindent\ignorespaces}
\def\applab#1{\DefWarn#1%
\xdef #1{\noexpand\hyperref{}{appendix}{\secn@m}{\secn@m}}%
\writedef{#1\leftbracket#1}\wrlabeL{#1=#1}}
\def\appsublab#1{\DefWarn#1%
\xdef #1{\noexpand\hyperref{}{appendix}{\secn@m.\the\subsecno}{\secn@m.\the\subsecno}}%
\writedef{#1\leftbracket#1}\wrlabeL{#1=#1}}
\newwrite\tfile \def\writetoca#1{}
\def\leaderfill{\leaders\hbox to 1em{\hss.\hss}\hfill}
% use this to write file with table of contents
\def\writetoc{\immediate\openout\tfile=\jobname.toc
   \def\writetoca##1{{\edef\next{\write\tfile{\noindent ##1
   \string\leaderfill{
% comment this line if you don't want hyperlinked page numbers on TOC
   \string\hyperref{}{page}{\noexpand\number\pageno}%
   {\noexpand\number\pageno}} \par}}\next}}
}
% and this lists table of contents on second pass
\newread\ch@ckfile
\def\listtoc{\immediate\closeout\tfile\immediate\openin\ch@ckfile=\jobname.toc
\ifeof\ch@ckfile\message{no file \jobname.toc, no table of contents this pass}%
\else\closein\ch@ckfile\centerline{\bf Contents}\nobreak\medskip%
{\baselineskip=16pt\footnotefont\parskip=0pt\catcode`\@=11\input\jobname.toc
\catcode`\@=12\bigbreak\bigskip}\fi}
\catcode`\@=12 % at signs are no longer letters
\font\ninerm=cmr9 \font\sixrm=cmr6 \font\ninei=cmmi9 \font\sixi=cmmi6
\font\ninesy=cmsy9 \font\sixsy=cmsy6 \font\ninebf=cmbx9
\font\nineit=cmti9 \font\ninesl=cmsl9 \skewchar\ninei='177
\skewchar\sixi='177 \skewchar\ninesy='60 \skewchar\sixsy='60
\def\ninepoint{\def\rm{\fam0\ninerm}% switch to footnote font
\textfont0=\ninerm \scriptfont0=\sixrm \scriptscriptfont0=\fiverm
\textfont1=\ninei \scriptfont1=\sixi \scriptscriptfont1=\fivei
\textfont2=\ninesy \scriptfont2=\sixsy \scriptscriptfont2=\fivesy
\textfont\itfam=\ninei \def\it{\fam\itfam\nineit}\def\sl{\fam\slfam\ninesl}%
\textfont\bffam=\ninebf \def\bf{\fam\bffam\ninebf}\rm}
%
%---------------------------------------------------------------------
\hyphenation{anom-aly anom-alies coun-ter-term coun-ter-terms}

%%%%%%%%%%%%%%% Subsubsection %%%%%%%%%%%%%%%%%%%%%%%%%%%%%%%%%%%%
\global\newcount\subsubsecno \global\subsubsecno=0
\def\subsubsec#1\par{\global\advance\subsubsecno by1%
{\toks0{#1}\message{(\the\secno\the\subsecno\the\subsubsecno. \the\toks0)}}%
\ifnum\lastpenalty>9000\else\bigbreak\fi
\noindent{\it\hyperdef\hypernoname{subsubsection}{\the\secno.\the\subsecno\the\subsubsecno}%
{\the\secno.\the\subsecno.\the\subsubsecno.} #1}
%%% Add Subsubsections to Index
%\writetoca{\string\quad{\string\hyperref{}{subsubsection}{\the\secno\the\subsecno\the
%\subsubsecno}{\baselineskip=9pt\it\the\secno.\the\subsecno.\the\subsubsecno.}}
% {\baselineskip=9pt\it\ #1}}
\par\nobreak\medskip\nobreak\noindent\ignorespaces}

% Caption for inline tikzpictures
\def\DefWarn#1{}
\def\tikzcaption#1#2{\DefWarn#1\xdef#1{Fig.~\the\figno}
\writedef{#1\leftbracket Fig.\noexpand~\the\figno}%
{
\smallskip
\leftskip=20pt \rightskip=20pt \baselineskip12pt\noindent
{{\bf Fig.~\the\figno}\ \ninepoint #2}
\bigskip
\global\advance\figno by1 \par}}

% convert numbers [1-9] to upper case letters [A-I]
\def\ntoalpha#1{%
\ifcase#1%
@%
\or A\or B\or C\or D\or E\or F\or G\or H\or I
\fi
}

% Appendix label
\global\newcount\appno \global\appno=1
\def\applab#1{\xdef #1{\ntoalpha\appno}\writedef{#1\leftbracket#1}\wrlabeL{#1=#1}
\global\advance\appno by1}

% Clean up the title page definitions
%\def\preprint#1 #2\par{\rightline{\vbox{\baselineskip12pt\hbox{#1}\hbox{#2}}}\vskip2cm}
% title with more than one line (note the blanck line in between)
%\title some line
%
%\tile another line
\def\title#1\par{\centerline{\bf #1}}
\def\author#1\par{\bigskip\bigskip\centerline{#1}}

\newcount\addressno

\def\email#1#2{%\unskip$^#1$
\footnote{\null}{\kern-\parindent \llap{$^#1$\hskip1pt}email: #2}}

% centermode for address lines
\def\startcenter{%
  \par
  \begingroup
  \leftskip=0pt plus 1fil
  \rightskip=\leftskip
  \parindent=0pt
  \parfillskip=0pt
}
\def\stopcenter{\endgroup}

\def\address{\bigskip%
  \ifnum\the\addressno=0\else\stopcenter\endgroup\fi
  \advance\addressno by 1%
  \begingroup
  \startcenter
  \it
  \obeylines
  \addressAux
}
\def\addressAux#1{#1}

% need to stop center mode and obeylines from address
\def\abstract{\stopcenter\endgroup\bigskip\bigskip\noindent}

% some sample definitions
\def\Dsl{\,\raise.15ex\hbox{/}\mkern-13.5mu D} %this one can be subscripted
\def\dsl{\raise.15ex\hbox{/}\kern-.57em\partial}
 
\def\boxeqn#1{\vcenter{\vbox{\hrule\hbox{\vrule\kern3pt\vbox{\kern3pt
	\hbox{${\displaystyle #1}$}\kern3pt}\kern3pt\vrule}\hrule}}}

 %pound sterling

\def\d{{\delta}}

\def\s{{\sigma}}

\def\half{{1\over 2}}

\def\({\left(}
\def\){\right)}

% blackboard bold

% primed summation symbol

% length of words, |P|

 %redefine plain TeX \Im..
% small inlined fractions, from the TeXbook
\def\sfrac#1/#2{\kern.1em\raise.5ex\hbox{\the\scriptfont0 #1}%
\kern-.1em/\kern-.15em\lower.25ex\hbox{\the\scriptfont0 #2}}

%shuffle product
\font\tenshuffle=shuffle10 \font\sevenshuffle=shuffle7 \font\fiveshuffle=shuffle7 at 5pt
\def\shuffle{{%
\def\Dshuffle{\mathbin{\hbox{\tenshuffle\char'001}}}%
\def\Sshuffle{\mathbin{\hbox{\sevenshuffle\char'001}}}%
\def\SSshuffle{\mathbin{\hbox{\fiveshuffle\char'001}}}%
\mathchoice{\Dshuffle}{\Dshuffle}{\Sshuffle}{\SSshuffle}}}

%\owedge

% From Knuth's \pfbox macro
\def\qed{\hbox{\hskip 3pt
%\lower2pt
\vbox{\hrule\hbox to 7pt{\vrule height 7pt\hfill\vrule}
\hrule}}\hskip3pt}

\def\PBT{{\rm PBT}}

% do not display overfull marks
\overfullrule=0pt\relax

\frenchspacing

% define labels in advance
\newread\instream \openin\instream= label.defs
\ifeof\instream \message{No labels in advance yet. Wait till next pass.}
\else \closein\instream \input label.defs
\fi
\writedefs

%%% References with hyperlinks to arxiv.org; both styles accepted
% Change arXiv to \arXiv ie
% [arXiv:hep-th/1234567].     --> [\arXiv:hep-th/1234567].
% [arXiv:1234.5678 [hep-th]]. --> [\arXiv:1234.5678 [hep-th]].
% Need to strip trailing [hep-th] (if present) to define valid URL
\def\arXiv:#1].{\hepthStrip#1 \nil}
\def\hepthStrip#1 #2\nil{\href{http://arxiv.org/abs/#1}{arXiv:#1 #2\unskip}].}

% from the paper 0102053 by Loday, to draw inlined PBTs
\catcode`@=11
\def\@height{height}
\def\@depth{depth}
\def\@width{width}

\newcount\@tempcnta
\newcount\@tempcntb

\newdimen\@tempdima
\newdimen\@tempdimb

\newbox\@tempboxa

\def\@ifnextchar#1#2#3{\let\@tempe #1\def\@tempa{#2}\def\@tempb{#3}\futurelet
    \@tempc\@ifnch}
\def\@ifnch{\ifx \@tempc \@sptoken \let\@tempd\@xifnch
      \else \ifx \@tempc \@tempe\let\@tempd\@tempa\else\let\@tempd\@tempb\fi
      \fi \@tempd}
\def\@ifstar#1#2{\@ifnextchar *{\def\@tempa*{#1}\@tempa}{#2}}

\def\@whilenoop#1{}
\def\@whilenum#1\do #2{\ifnum #1\relax #2\relax\@iwhilenum{#1\relax 
     #2\relax}\fi}
\def\@iwhilenum#1{\ifnum #1\let\@nextwhile=\@iwhilenum 
         \else\let\@nextwhile=\@whilenoop\fi\@nextwhile{#1}}

\def\@whiledim#1\do #2{\ifdim #1\relax#2\@iwhiledim{#1\relax#2}\fi}
\def\@iwhiledim#1{\ifdim #1\let\@nextwhile=\@iwhiledim 
        \else\let\@nextwhile=\@whilenoop\fi\@nextwhile{#1}}

\newdimen\@wholewidth
\newdimen\@halfwidth
\newdimen\unitlength \unitlength =1pt
\newbox\@picbox
\newdimen\@picht

\def\@nnil{\@nil}
\def\@empty{}
\def\@fornoop#1\@@#2#3{}

\def\@for#1:=#2\do#3{\edef\@fortmp{#2}\ifx\@fortmp\@empty \else
    \expandafter\@forloop#2,\@nil,\@nil\@@#1{#3}\fi}

\def\@forloop#1,#2,#3\@@#4#5{\def#4{#1}\ifx #4\@nnil \else
       #5\def#4{#2}\ifx #4\@nnil \else#5\@iforloop #3\@@#4{#5}\fi\fi}

\def\@iforloop#1,#2\@@#3#4{\def#3{#1}\ifx #3\@nnil 
       \let\@nextwhile=\@fornoop \else
      #4\relax\let\@nextwhile=\@iforloop\fi\@nextwhile#2\@@#3{#4}}

\def\@tfor#1:=#2\do#3{\xdef\@fortmp{#2}\ifx\@fortmp\@empty \else
    \@tforloop#2\@nil\@nil\@@#1{#3}\fi}
\def\@tforloop#1#2\@@#3#4{\def#3{#1}\ifx #3\@nnil 
       \let\@nextwhile=\@fornoop \else
      #4\relax\let\@nextwhile=\@tforloop\fi\@nextwhile#2\@@#3{#4}}

\def\@makepicbox(#1,#2){\leavevmode\@ifnextchar 
   [{\@imakepicbox(#1,#2)}{\@imakepicbox(#1,#2)[]}}

\long\def\@imakepicbox(#1,#2)[#3]#4{\vbox to#2\unitlength
   {\let\mb@b\vss \let\mb@l\hss\let\mb@r\hss
    \let\mb@t\vss
    \@tfor\@tempa :=#3\do{\expandafter\let
        \csname mb@\@tempa\endcsname\relax}%
\mb@t\hbox to #1\unitlength{\mb@l #4\mb@r}\mb@b}}

\def\picture(#1,#2){\@ifnextchar({\@picture(#1,#2)}{\@picture(#1,#2)(0,0)}}

\def\@picture(#1,#2)(#3,#4){\@picht #2\unitlength
\setbox\@picbox\hbox to #1\unitlength\bgroup 
\hskip -#3\unitlength \lower #4\unitlength \hbox\bgroup\ignorespaces}

\def\endpicture{\egroup\hss\egroup\ht\@picbox\@picht
\dp\@picbox\z@\leavevmode\box\@picbox}

\long\def\put(#1,#2)#3{\@killglue\raise#2\unitlength\hbox to \z@{\kern
#1\unitlength #3\hss}\ignorespaces}

\long\def\multiput(#1,#2)(#3,#4)#5#6{\@killglue\@multicnt=#5\relax
\@xdim=#1\unitlength
\@ydim=#2\unitlength
\@whilenum \@multicnt > 0\do
{\raise\@ydim\hbox to \z@{\kern
\@xdim #6\hss}\advance\@multicnt \m@ne\advance\@xdim
#3\unitlength\advance\@ydim #4\unitlength}\ignorespaces}

\def\@killglue{\unskip\@whiledim \lastskip >\z@\do{\unskip}}

\def\thinlines{\let\@linefnt\tenln \let\@circlefnt\tencirc
  \@wholewidth\fontdimen8\tenln \@halfwidth .5\@wholewidth}
\def\thicklines{\let\@linefnt\tenlnw \let\@circlefnt\tencircw
  \@wholewidth\fontdimen8\tenlnw \@halfwidth .5\@wholewidth}

\def\linethickness#1{\@wholewidth #1\relax \@halfwidth .5\@wholewidth}

\def\shortstack{\@ifnextchar[{\@shortstack}{\@shortstack[c]}}

\def\@shortstack[#1]{\leavevmode
\vbox\bgroup\baselineskip-1pt\lineskip 3pt\let\mb@l\hss
\let\mb@r\hss \expandafter\let\csname mb@#1\endcsname\relax
\let\\\@stackcr\@ishortstack}

\def\@ishortstack#1{\halign{\mb@l ##\unskip\mb@r\cr #1\crcr}\egroup}

\def\@stackcr{\@ifstar{\@ixstackcr}{\@ixstackcr}}
\def\@ixstackcr{\@ifnextchar[{\@istackcr}{\cr\ignorespaces}}

\def\@istackcr[#1]{\cr\noalign{\vskip #1}\ignorespaces}

\newif\if@negarg

\def\droite(#1,#2)#3{\@xarg #1\relax \@yarg #2\relax
\@linelen=#3\unitlength
\ifnum\@xarg =0 \@vline 
  \else \ifnum\@yarg =0 \@hline \else \@sline\fi
\fi}

\def\@sline{\ifnum\@xarg< 0 \@negargtrue \@xarg -\@xarg \@yyarg -\@yarg
  \else \@negargfalse \@yyarg \@yarg \fi
\ifnum \@yyarg >0 \@tempcnta\@yyarg \else \@tempcnta -\@yyarg \fi
\ifnum\@tempcnta>6 \@badlinearg\@tempcnta0 \fi
\ifnum\@xarg>6 \@badlinearg\@xarg 1 \fi
\setbox\@linechar\hbox{\@linefnt\@getlinechar(\@xarg,\@yyarg)}%
\ifnum \@yarg >0 \let\@upordown\raise \@clnht\z@
   \else\let\@upordown\lower \@clnht \ht\@linechar\fi
\@clnwd=\wd\@linechar
\if@negarg \hskip -\wd\@linechar \def\@tempa{\hskip -2\wd\@linechar}\else
     \let\@tempa\relax \fi
\@whiledim \@clnwd <\@linelen \do
  {\@upordown\@clnht\copy\@linechar
   \@tempa
   \advance\@clnht \ht\@linechar
   \advance\@clnwd \wd\@linechar}%
\advance\@clnht -\ht\@linechar
\advance\@clnwd -\wd\@linechar
\@tempdima\@linelen\advance\@tempdima -\@clnwd
\@tempdimb\@tempdima\advance\@tempdimb -\wd\@linechar
\if@negarg \hskip -\@tempdimb \else \hskip \@tempdimb \fi
\multiply\@tempdima \@m
\@tempcnta \@tempdima \@tempdima \wd\@linechar \divide\@tempcnta \@tempdima
\@tempdima \ht\@linechar \multiply\@tempdima \@tempcnta
\divide\@tempdima \@m
\advance\@clnht \@tempdima
\ifdim \@linelen <\wd\@linechar
   \hskip \wd\@linechar
  \else\@upordown\@clnht\copy\@linechar\fi}

\def\@hline{\ifnum \@xarg <0 \hskip -\@linelen \fi
\vrule \@height \@halfwidth \@depth \@halfwidth \@width \@linelen
\ifnum \@xarg <0 \hskip -\@linelen \fi}

\def\@getlinechar(#1,#2){\@tempcnta#1\relax\multiply\@tempcnta 8
\advance\@tempcnta -9 \ifnum #2>0 \advance\@tempcnta #2\relax\else
\advance\@tempcnta -#2\relax\advance\@tempcnta 64 \fi
\char\@tempcnta}

\def\vector(#1,#2)#3{\@xarg #1\relax \@yarg #2\relax
\@tempcnta \ifnum\@xarg<0 -\@xarg\else\@xarg\fi
\ifnum\@tempcnta<5\relax
\@linelen=#3\unitlength
\ifnum\@xarg =0 \@vvector 
  \else \ifnum\@yarg =0 \@hvector \else \@svector\fi
\fi
\else\@badlinearg\fi}

\def\@hvector{\@hline\hbox to 0pt{\@linefnt 
\ifnum \@xarg <0 \@getlarrow(1,0)\hss\else
    \hss\@getrarrow(1,0)\fi}}

\def\@vvector{\ifnum \@yarg <0 \@downvector \else \@upvector \fi}

\def\@svector{\@sline
\@tempcnta\@yarg \ifnum\@tempcnta <0 \@tempcnta=-\@tempcnta\fi
\ifnum\@tempcnta <5
  \hskip -\wd\@linechar
  \@upordown\@clnht \hbox{\@linefnt  \if@negarg 
  \@getlarrow(\@xarg,\@yyarg) \else \@getrarrow(\@xarg,\@yyarg) \fi}%
\else\@badlinearg\fi}

\def\@getlarrow(#1,#2){\ifnum #2 =\z@ \@tempcnta='33\else
\@tempcnta=#1\relax\multiply\@tempcnta \sixt@@n \advance\@tempcnta
-9 \@tempcntb=#2\relax\multiply\@tempcntb \tw@
\ifnum \@tempcntb >0 \advance\@tempcnta \@tempcntb\relax
\else\advance\@tempcnta -\@tempcntb\advance\@tempcnta 64
\fi\fi\char\@tempcnta}

\def\@getrarrow(#1,#2){\@tempcntb=#2\relax
\ifnum\@tempcntb < 0 \@tempcntb=-\@tempcntb\relax\fi
\ifcase \@tempcntb\relax \@tempcnta='55 \or 
\ifnum #1<3 \@tempcnta=#1\relax\multiply\@tempcnta
24 \advance\@tempcnta -6 \else \ifnum #1=3 \@tempcnta=49
\else\@tempcnta=58 \fi\fi\or 
\ifnum #1<3 \@tempcnta=#1\relax\multiply\@tempcnta
24 \advance\@tempcnta -3 \else \@tempcnta=51\fi\or 
\@tempcnta=#1\relax\multiply\@tempcnta
\sixt@@n \advance\@tempcnta -\tw@ \else
\@tempcnta=#1\relax\multiply\@tempcnta
\sixt@@n \advance\@tempcnta 7 \fi\ifnum #2<0 \advance\@tempcnta 64 \fi
\char\@tempcnta}

\def\@vline{\ifnum \@yarg <0 \@downline \else \@upline\fi}

\def\@upline{\hbox to \z@{\hskip -\@halfwidth \vrule \@width \@wholewidth
   \@height \@linelen \@depth \z@\hss}}

\def\@downline{\hbox to \z@{\hskip -\@halfwidth \vrule \@width \@wholewidth
   \@height \z@ \@depth \@linelen \hss}}

\def\@upvector{\@upline\setbox\@tempboxa\hbox{\@linefnt\char'66}\raise 
     \@linelen \hbox to\z@{\lower \ht\@tempboxa\box\@tempboxa\hss}}

\def\@downvector{\@downline\lower \@linelen
      \hbox to \z@{\@linefnt\char'77\hss}}

\def\dashbox#1(#2,#3){\leavevmode\hbox to \z@{\baselineskip \z@%
\lineskip \z@%
\@dashdim=#2\unitlength%
\@dashcnt=\@dashdim \advance\@dashcnt 200
\@dashdim=#1\unitlength\divide\@dashcnt \@dashdim
\ifodd\@dashcnt\@dashdim=\z@%
\advance\@dashcnt \@ne \divide\@dashcnt \tw@ 
\else \divide\@dashdim \tw@ \divide\@dashcnt \tw@
\advance\@dashcnt \m@ne
\setbox\@dashbox=\hbox{\vrule \@height \@halfwidth \@depth \@halfwidth
\@width \@dashdim}\put(0,0){\copy\@dashbox}%
\put(0,#3){\copy\@dashbox}%
\put(#2,0){\hskip-\@dashdim\copy\@dashbox}%
\put(#2,#3){\hskip-\@dashdim\box\@dashbox}%
\multiply\@dashdim 3 
\fi
\setbox\@dashbox=\hbox{\vrule \@height \@halfwidth \@depth \@halfwidth
\@width #1\unitlength\hskip #1\unitlength}\@tempcnta=0
\put(0,0){\hskip\@dashdim \@whilenum \@tempcnta <\@dashcnt
\do{\copy\@dashbox\advance\@tempcnta \@ne }}\@tempcnta=0
\put(0,#3){\hskip\@dashdim \@whilenum \@tempcnta <\@dashcnt
\do{\copy\@dashbox\advance\@tempcnta \@ne }}%
\@dashdim=#3\unitlength%
\@dashcnt=\@dashdim \advance\@dashcnt 200
\@dashdim=#1\unitlength\divide\@dashcnt \@dashdim
\ifodd\@dashcnt \@dashdim=\z@%
\advance\@dashcnt \@ne \divide\@dashcnt \tw@
\else
\divide\@dashdim \tw@ \divide\@dashcnt \tw@
\advance\@dashcnt \m@ne
\setbox\@dashbox\hbox{\hskip -\@halfwidth
\vrule \@width \@wholewidth 
\@height \@dashdim}\put(0,0){\copy\@dashbox}%
\put(#2,0){\copy\@dashbox}%
\put(0,#3){\lower\@dashdim\copy\@dashbox}%
\put(#2,#3){\lower\@dashdim\copy\@dashbox}%
\multiply\@dashdim 3
\fi
\setbox\@dashbox\hbox{\vrule \@width \@wholewidth 
\@height #1\unitlength}\@tempcnta0
\put(0,0){\hskip -\@halfwidth \vbox{\@whilenum \@tempcnta < \@dashcnt
\do{\vskip #1\unitlength\copy\@dashbox\advance\@tempcnta \@ne }%
\vskip\@dashdim}}\@tempcnta0
\put(#2,0){\hskip -\@halfwidth \vbox{\@whilenum \@tempcnta< \@dashcnt
\relax\do{\vskip #1\unitlength\copy\@dashbox\advance\@tempcnta \@ne }%
\vskip\@dashdim}}}\@makepicbox(#2,#3)}

\newif\if@ovt 
\newif\if@ovb 
\newif\if@ovl 
\newif\if@ovr 
\newdimen\@ovxx
\newdimen\@ovyy
\newdimen\@ovdx
\newdimen\@ovdy
\newdimen\@ovro
\newdimen\@ovri

\def\@getcirc#1{\@tempdima #1\relax \advance\@tempdima 2pt\relax
  \@tempcnta\@tempdima
  \@tempdima 4pt\relax \divide\@tempcnta\@tempdima
  \ifnum \@tempcnta > 10\relax \@tempcnta 10\relax\fi
  \ifnum \@tempcnta >\z@ \advance\@tempcnta\m@ne
    \else \@warning{Oval too small}\fi
  \multiply\@tempcnta 4\relax
  \setbox \@tempboxa \hbox{\@circlefnt
  \char \@tempcnta}\@tempdima \wd \@tempboxa}

\def\@put#1#2#3{\raise #2\hbox to \z@{\hskip #1#3\hss}}

\def\oval(#1,#2){\@ifnextchar[{\@oval(#1,#2)}{\@oval(#1,#2)[]}}

\def\@oval(#1,#2)[#3]{\begingroup\boxmaxdepth \maxdimen
  \@ovttrue \@ovbtrue \@ovltrue \@ovrtrue
  \@tfor\@tempa :=#3\do{\csname @ov\@tempa false\endcsname}\@ovxx
  #1\unitlength \@ovyy #2\unitlength
  \@tempdimb \ifdim \@ovyy >\@ovxx \@ovxx\else \@ovyy \fi
  \advance \@tempdimb -2pt\relax  %%%% added 7 Dec 89
  \@getcirc \@tempdimb
  \@ovro \ht\@tempboxa \@ovri \dp\@tempboxa
  \@ovdx\@ovxx \advance\@ovdx -\@tempdima \divide\@ovdx \tw@
  \@ovdy\@ovyy \advance\@ovdy -\@tempdima \divide\@ovdy \tw@
  \@circlefnt \setbox\@tempboxa
  \hbox{\if@ovr \@ovvert32\kern -\@tempdima \fi
  \if@ovl \kern \@ovxx \@ovvert01\kern -\@tempdima \kern -\@ovxx \fi
  \if@ovt \@ovhorz \kern -\@ovxx \fi
  \if@ovb \raise \@ovyy \@ovhorz \fi}\advance\@ovdx\@ovro
  \advance\@ovdy\@ovro \ht\@tempboxa\z@ \dp\@tempboxa\z@
  \@put{-\@ovdx}{-\@ovdy}{\box\@tempboxa}%
  \endgroup}

\def\@ovvert#1#2{\vbox to \@ovyy{%
    \if@ovb \@tempcntb \@tempcnta \advance \@tempcntb by #1\relax
      \kern -\@ovro \hbox{\char \@tempcntb}\nointerlineskip
    \else \kern \@ovri \kern \@ovdy \fi
    \leaders\vrule width \@wholewidth\vfil \nointerlineskip
    \if@ovt \@tempcntb \@tempcnta \advance \@tempcntb by #2\relax
      \hbox{\char \@tempcntb}%
    \else \kern \@ovdy \kern \@ovro \fi}}

\def\@ovhorz{\hbox to \@ovxx{\kern \@ovro
    \if@ovr \else \kern \@ovdx \fi
    \leaders \hrule height \@wholewidth \hfil
    \if@ovl \else \kern \@ovdx \fi
    \kern \@ovri}}

\def\circle{\@ifstar{\@dot}{\@circle}}
\def\@circle#1{\begingroup \boxmaxdepth \maxdimen \@tempdimb #1\unitlength
   \ifdim \@tempdimb >15.5pt\relax \@getcirc\@tempdimb
      \@ovro\ht\@tempboxa 
     \setbox\@tempboxa\hbox{\@circlefnt
      \advance\@tempcnta\tw@ \char \@tempcnta
      \advance\@tempcnta\m@ne \char \@tempcnta \kern -2\@tempdima
      \advance\@tempcnta\tw@
      \raise \@tempdima \hbox{\char\@tempcnta}\raise \@tempdima
        \box\@tempboxa}\ht\@tempboxa\z@ \dp\@tempboxa\z@
      \@put{-\@ovro}{-\@ovro}{\box\@tempboxa}%
   \else  \@circ\@tempdimb{96}\fi\endgroup}

\def\@dot#1{\@tempdimb #1\unitlength \@circ\@tempdimb{112}}

\def\@circ#1#2{\@tempdima #1\relax \advance\@tempdima .5pt\relax
   \@tempcnta\@tempdima \@tempdima 1pt\relax
   \divide\@tempcnta\@tempdima 
   \ifnum\@tempcnta > 15\relax \@tempcnta 15\relax \fi    
   \ifnum \@tempcnta >\z@ \advance\@tempcnta\m@ne\fi
   \advance\@tempcnta #2\relax
   \@circlefnt \char\@tempcnta}

%INITIALIZATION
\font\tenln line10
%\font\tencirc circle10
\font\tencirc lcircle10
\font\tenlnw linew10
%\font\tencircw circlew10
\font\tencircw lcirclew10

\thinlines   

\newcount\@xarg
\newcount\@yarg
\newcount\@yyarg
\newcount\@multicnt 
\newdimen\@xdim
\newdimen\@ydim
\newbox\@linechar
\newdimen\@linelen
\newdimen\@clnwd
\newdimen\@clnht
\newdimen\@dashdim
\newbox\@dashbox
\newcount\@dashcnt
\catcode`@=12
%%% Local Variables: 
%%% mode: plain-tex
%%% TeX-master: t
%%% End: 

%\input poorPict.tex
\catcode`@=11
\font\@linefnt linew10 at 2.4pt
\catcode`@=12

\def\arbreun{\kern-0.4ex
\hbox{\unitlength=.25pt
\picture(60,40)(0,0)
\put(30,0){\droite(0,1){20}}
\put(30,20){\droite(-1,1){30}}
\put(30,20){\droite(1,1){30}}
\put(20,30){\droite(1,1){20}}
\put(10,40){\droite(1,1){10}}
\endpicture}\kern 0.4ex}

\def\arbredeux{\kern-0.4ex
\hbox{\unitlength=.25pt
\picture(60,40)(0,0)
\put(30,0){\droite(0,1){20}}
\put(30,20){\droite(-1,1){30}}
\put(30,20){\droite(1,1){30}}
\put(20,30){\droite(1,1){20}}
\put(30,40){\droite(-1,1){10}}
\endpicture}\kern 0.4ex}

\def\arbretrois{\kern-0.4ex
\hbox{\unitlength=.25pt
\picture(60,40)(0,0)
\put(30,0){\droite(0,1){20}}
\put(30,20){\droite(-1,1){30}}
\put(30,20){\droite(1,1){30}}
\put(50,40){\droite(-1,1){10}}
\put(10,40){\droite(1,1){10}}
\endpicture}\kern 0.4ex}

\def\arbrequatre{\kern-0.4ex
\hbox{\unitlength=.25pt
\picture(60,40)(0,0)
\put(30,0){\droite(0,1){20}}
\put(30,20){\droite(-1,1){30}}
\put(30,20){\droite(1,1){30}}
\put(40,30){\droite(-1,1){20}}
\put(30,40){\droite(1,1){10}}
\endpicture}\kern 0.4ex}

\def\arbrecinq{\kern-0.4ex
\hbox{\unitlength=.25pt
\picture(60,40)(0,0)
\put(30,0){\droite(0,1){20}}
\put(30,20){\droite(-1,1){30}}
\put(30,20){\droite(1,1){30}}
\put(40,30){\droite(-1,1){20}}
\put(50,40){\droite(-1,1){10}}
\endpicture}\kern 0.4ex}

\def\arbreA{\kern-0.4ex
\hbox{\unitlength=.25pt
\picture(60,40)(0,0)
\put(30,0){\droite(0,1){20}}
\put(30,20){\droite(-1,1){30}}
\put(30,20){\droite(1,1){30}}
\endpicture}\kern 0.4ex}

\def\arbreB{\kern-0.4ex
\hbox{\unitlength=.25pt
\picture(60,40)(0,0)
\put(30,0){\droite(0,1){20}}
\put(30,20){\droite(-1,1){30}}
\put(30,20){\droite(1,1){30}}
\put(15,35){\droite(1,1){15}}
\endpicture}\kern 0.4ex}

\def\arbreC{\kern-0.4ex
\hbox{\unitlength=.25pt
\picture(60,40)(0,0)
\put(30,0){\droite(0,1){20}}
\put(30,20){\droite(-1,1){30}}
\put(30,20){\droite(1,1){30}}
\put(45,35){\droite(-1,1){15}}
\endpicture}\kern 0.4ex}

\lref\mouldec{
	J. Ecalle. ``ARI/GARI, la dimorphie et l’arithmétique des
	multizêtas: un premier bilan''. J.
	Théor. Nombres Bordeaux, 15(2):411–478, 2003.
}

\lref\BerendsME{
	F.A.~Berends and W.T.~Giele,
  	``Recursive Calculations for Processes with n Gluons,''
	Nucl.\ Phys.\ B {\bf 306}, 759 (1988).
	%%CITATION = Print-88-0100 (LEIDEN)%%
}
\lref\BGPS{
	C.R.~Mafra and O.~Schlotterer,
  	``Berends-Giele recursions and the BCJ duality in superspace and components,''
	JHEP {\bf 1603}, 097 (2016).
	[arXiv:1510.08846 [hep-th]].
	%%CITATION = DAMTP-2015-69%%
}

\lref\patras{
	F.~Patras, C.~Reutenauer, M.~Schocker, ``On the Garsia Lie Idempotent'',
   	Canad. Math. Bull. {\bf 48} (2005), 445-454
}

\lref\barcelo{
	Barcelo, H. and Sundaram, S., ``On Some Submodules of the Action of the Symmetrical
	Group on the Free Lie Algebra''. Journal of Algebra, 154(1), (1993) pp.12-26.
}

\lref\NLSM{
	J.~J.~M.~Carrasco, C.R.~Mafra and O.~Schlotterer,
  	``Abelian Z-theory: NLSM amplitudes and $\alpha$'-corrections from the open string,''
	JHEP {\bf 1706}, 093 (2017).
	[arXiv:1608.02569 [hep-th]].
	%%CITATION = arXiv:1608.02569%%
}

\lref\psf{
 	N.~Berkovits,
	``Super-Poincare covariant quantization of the superstring,''
	JHEP {\bf 0004}, 018 (2000)
	[arXiv:hep-th/0001035].
	%%CITATION = JHEPA,0004,018;%%
}

\lref\Polylogs{
	J.~Broedel, O.~Schlotterer and S.~Stieberger,
	``Polylogarithms, Multiple Zeta Values and Superstring Amplitudes,''
	Fortsch.\ Phys.\  {\bf 61}, 812 (2013).
	[arXiv:1304.7267 [hep-th]].
	%%CITATION = DAMTP-2013-22%%
}

\lref\BGschocker{
	M. Schocker,
	``Lie elements and Knuth relations,'' Canad. J. Math. {\bf 56} (2004), 871-882.
	[math/0209327].
}

\lref\Gauge{
	S.~Lee, C.R.~Mafra and O.~Schlotterer,
  	``Non-linear gauge transformations in $D=10$ SYM theory and the BCJ duality,''
	JHEP {\bf 1603}, 090 (2016).
	[arXiv:1510.08843 [hep-th]].
	%%CITATION = DAMTP-2015-68%%
}
\lref\FTlimit{
	C.R.~Mafra,
  	``Berends-Giele recursion for double-color-ordered amplitudes,''
	JHEP {\bf 1607}, 080 (2016).
	[arXiv:1603.09731 [hep-th]].
	%%CITATION = arXiv:1603.09731%%
}
\lref\KKref{
	R.~Kleiss and H.~Kuijf,
	``Multi - Gluon Cross-sections and Five Jet Production at Hadron Colliders,''
	Nucl.\ Phys.\ B {\bf 312}, 616 (1989)..
	%%CITATION = Print-88-0425 (LEIDEN)%%
}
\lref\KLTref{
	H.~Kawai, D.~C.~Lewellen and S.~H.~H.~Tye,
  	``A Relation Between Tree Amplitudes of Closed and Open Strings,''
	Nucl.\ Phys.\ B {\bf 269}, 1 (1986)..
	%%CITATION = CLNS-85/667%%
}
\lref\oldMomKer{
	Z.~Bern, L.~J.~Dixon, M.~Perelstein and J.~S.~Rozowsky,
	``Multileg one loop gravity amplitudes from gauge theory,''
	Nucl.\ Phys.\ B {\bf 546}, 423 (1999).
	[hep-th/9811140].
	%%CITATION = hep-th/9811140%%
}
\lref\MomKer{
	N.~E.~J.~Bjerrum-Bohr, P.~H.~Damgaard, T.~Sondergaard and P.~Vanhove,
	``The Momentum Kernel of Gauge and Gravity Theories,''
	JHEP {\bf 1101}, 001 (2011).
	[arXiv:1010.3933 [hep-th]].
	%%CITATION = arXiv:1010.3933%%
}
\lref\nptMethod{
	C.R.~Mafra, O.~Schlotterer, S.~Stieberger and D.~Tsimpis,
	``A recursive method for SYM n-point tree amplitudes,''
	Phys.\ Rev.\ D {\bf 83}, 126012 (2011).
	[arXiv:1012.3981 [hep-th]].
	%%CITATION = arXiv:1012.3981%%
}
\lref\loday{
	J.-L~Loday and M.~O.~Ronco, ``Hopf algebra of the planar binary trees'',
	Adv. Math. 139 (1998), no. 2, 293--309.
}
\lref\vallette{
	J.-L. Loday and B. Vallette, ``Algebraic operads'',
	Grundlehren Math. Wiss. 346, Springer, Heidelberg, 2012.
}
\lref\stanley{
	R.P. Stanley, ``Enumerative Combinatorics'', vols. I and II, second edition, Cambridge, UK: Univ. Pr.
	(2012)
}
\lref\reutenauer{
	C.~Reutenauer,
	``Free Lie Algebras'', London Mathematical Society Monographs, 1993.
}
\lref\EOMBBs{
	C.R.~Mafra and O.~Schlotterer,
  	``Multiparticle SYM equations of motion and pure spinor BRST blocks,''
	JHEP {\bf 1407}, 153 (2014).
	[arXiv:1404.4986 [hep-th]].
	%%CITATION = AEI-2014-011%%
}
\lref\Ree{
	R. Ree, ``Lie elements and an algebra associated with shuffles'', Ann.
	Math. {\bf 62}, No. 2 (1958),
	210--220.
}
\lref\chapoton{
	F.~Chapoton F. ``The anticyclic operad of moulds''. International Mathematics
	Research Notices. 2007 Jan 1;2007, math/0609436
}
\lref\garsia{
	A.M. Garsia, ``Combinatorics of the Free Lie Algebra and the Symmetric Group'',
	In Analysis, et Cetera, edited by Paul H. Rabinowitz and Eduard Zehnder,
	Academic Press, (1990) 309-382
}
\lref\DPellis{
	F.~Cachazo, S.~He and E.Y.~Yuan,
	``Scattering of Massless Particles: Scalars, Gluons and Gravitons,''
	JHEP {\bf 1407}, 033 (2014).
	[arXiv:1309.0885 [hep-th]].
	%%CITATION = arXiv:1309.0885%%
}
\lref\cycfac{
	Moszkowski, Paul. ``A solution to a problem of D\'enes:
	a bijection between trees and factorizations of cyclic permutations.'' European
	Journal of Combinatorics 10, no. 1 (1989): 13-16.
	\semi
	Goulden, Ian P., and S. Pepper. ``Labelled trees and
	factorizations of a cycle into transpositions.''
	Discrete Mathematics 113, no. 1-3 (1993): 263-268.
}
\lref\nptString{
	C.R.~Mafra, O.~Schlotterer and S.~Stieberger,
  	``Complete N-Point Superstring Disk Amplitude I. Pure Spinor Computation,''
	Nucl.\ Phys.\ B {\bf 873}, 419 (2013).
	[arXiv:1106.2645 [hep-th]].
	%%CITATION = arXiv:1106.2645%%
}
\lref\dialg{
	Loday, J.L., ``Dialgebras''. In Dialgebras and related operads (2001) 7-66.
	Springer, Berlin, Heidelberg.
	[arXiv:math/0102053 [math.QA]].
}

\title Planar binary trees in scattering amplitudes

\author
Carlos R. Mafra\email{\dagger}{c.r.mafra@soton.ac.uk}

\address
Mathematical Sciences and STAG Research Centre, University of Southampton,
Highfield, Southampton, SO17 1BJ, UK

\abstract
These notes are a written version of my talk given at the CARMA workshop in June
2017, with some additional material. I presented a few concepts that have
recently been used in the computation of tree-level scattering amplitudes
(mostly using pure spinor methods but not restricted to it) in a context that
could be of interest to the combinatorics community. In particular, I focused on
the appearance of {\it planar binary trees} in scattering amplitudes and
presented some curious identities obeyed by related objects, some of which are
known to be true only via explicit examples.

%**************************
\newsec Planar binary trees

The basic ingredients in the following discussions are the {\it planar binary
trees} (pb trees). Recall that a planar tree is binary if every vertex is cubic
(or trivalent), with one {\it root} and two {\it leaves}. It is customary to
denote by ${\rm PBT}_n$ the set containing all planar binary trees with $n$
leaves. Their sizes are given by the Catalan numbers; $C_{n-1}=1,1,2,5,14,42,\ldots$
For example \dialg,
$$
\PBT_1 = \{ |\}\,,
\PBT_2= \{\  \arbreA \  \}\,,
\PBT_3= \{\ \arbreB ,\arbreC \  \}\,,
\PBT_4 = \{\ \arbreun,\arbredeux,\arbretrois ,\arbrequatre ,\arbrecinq \ \}\,.
$$
Let us now successively add more structure to planar binary trees.
The motivation for doing this comes from the physics of scattering amplitudes
but for the moment let us focus on their
intrinsic combinatorial value.

In the subsequent discussions {\it words} are composed of permutations
from the alphabet of
natural numbers and will be written in upper case (e.g. $A=14532$)
while {\it letters} will be written in lower case (e.g. $i=3$).
The length of a word $A$ is denoted $|A|$. The generalized momenta
$k_A^m$ and the generalized Mandelstam variables $s_P$ are defined
as
\eqn\defs{
k_A^m = k^m_{a_1} + k^m_{a_2} + \cdots + k^m_{a_{|A|}},\quad
s_{P}\equiv \half k_P\cdot k_P\,,\quad
k_i\cdot k_j \equiv s_{ij} = s_{ji},
}
where the momentum for a single letter squares to zero, $k_i\cdot k_i = 0$.
For example $k_{12}\cdot k_3 = (k_1\cdot k_3 + k_2\cdot k_3) = s_{13}+s_{23}$
as well as $s_{123} = s_{12} + s_{13} + s_{23}$.
In addition,
every labelled object or function will be considered linear in words, for
example $T_{123+321}\equiv T_{123}+T_{321}$. This is also extended to cases
such as
\eqn\funlin{
b\Big({123\over s_{12}} + {231\over s_{23}}\Big)\equiv
{1\over s_{12}}b(123) + {1\over s_{23}}b(231)\,.
}

%***************************************************
\subsec Planar binary trees and Mandelstam variables

Let us associate to each pb tree a rational function of $s_{ij}$ following
a recursive setup similar to that of Garsia \garsia:
to each word $P$ from
the alphabet $X=\{1,2, \ldots,n\}$ we let the pair $(P,T)$ represent the tree in
which the letters of $P$ are successively assigned to the leaves of $T$ from
left to right. Now let $T_1$ and $T_2$ be the left and right subtrees of $T$ and
let $P_1$ and $P_2$ be their corresponding words such that $P=P_1P_2$. The map
from the tree $T$ to a function of $s_{ij}$ is defined by,
\eqn\bPT{
\phi(P,T) \equiv {1\over s_P} \phi(P_1,T_1) \phi(P_2,T_2),\qquad \phi(i,T)\equiv i
}
For example, the action of the map \bPT\ on the two pb trees in $\PBT_3$ is
depicted in \figsingle.
\ifig\figsingle{Example applications of the map $\phi(P,T)$ defined in \bPT.}
{\epsfxsize=0.38\hsize\epsfbox{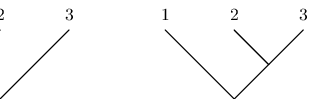}}

Now let me illustrate a common theme in the discussions to follow.
Suppose we are interested in the
image of the map \bPT\ not for an individual tree but for the sum
over all pb trees in $\PBT_n$. While it is straightforward to sum
$\phi(P,T_j)$ over all pb trees in $\PBT_n$, there is another way
to get the answer: drop
the tree specification $\phi(P,T)\equiv\phi(P)$ and
evaluate the following recursion \FTlimit,
\eqn\BGtree{
\phi(P) \equiv {1\over s_P}\!\!\sum_{XY=P}\!\!\phi(X)\phi(Y)\,,\qquad
\phi(i) \equiv i\,,\qquad\phi(\emptyset)\equiv 0\,,
}
where $\sum_{XY=P}$ denotes the sum over all deconcatenations of the word $P$
into $X$ and $Y$.
If we denote by
$C_n$ the number of terms in the expansion of $\phi(12 \ldots n{+}1)$,
it is easy to see that \BGtree\ gives rise to the
recurrence relation for the Catalan numbers,
$C_0=1$, $C_{n+1} =\sum_{i=0}^n C_i C_{n-i}$. Therefore all the
pb trees with $n$ leaves can be generated by considering the deconcatenation
of words of length $n$.
The first few expansions of \BGtree\ are given by
\eqnn\bintrees
$$\displaylines{
\phi(1) = 1,\quad\quad \phi(12) = {1\over s_{12}}, \quad \quad
\phi(123) = {1 \over s_{12} s_{123}} + { 1 \over s_{23}s_{123}}\,,\hfil\bintrees\hfilneg\cr
\phi(1234) = {1 \over s_{1234}} \Big( {1\over s_{12}s_{123} }
 + {1\over s_{23}s_{123} } + {1 \over s_{12}s_{34} }
 + {1 \over s_{34}s_{234} } + {1 \over s_{23}s_{234} } \Big)\,.
}$$
The definition \bPT\ resembles the map $\psi(T)$
defined by Chapoton in Proposition 3.2 of \chapoton, but
unfortunately they are not equivalent\foot{In \chapoton\ the pb trees
are mapped to variables $1/u_i$ defined by the {\it intervals} between sequential
leaves, while the Mandelstam variables $s_P$ may contain arbitrary leaves (such as
$s_{13}$).}
%consider the sum of the two expressions in $PBT_2$; while
%in \chapoton\ one gets $1/u_1u_2$ in terms of the variables $u_i$ labelled
%by the {\it intervals} between leaves, here we get
%$(1/s_{123}s_{12}+1/s_{23}s_{123})$ which does not simplify to a single term.}.
In any case, as discussed in \chapoton\ the map $\psi(T)$
gives rise to a {\it mould} \mouldec, so one may ask similar questions here.
As we will see below, one can modify the recursion \BGtree\ to obtain an
{\it alternal mould} of planar binary trees.

%**************************************************
\subsec An alternal mould of planar binary trees

A variation of the recursion \BGtree\ gives rise to
an {\it alternal mould of planar binary trees}.
Define $\phi(P|Q)$ in terms of {\it two} words
$P$ and $Q$ recursively as \FTlimit,
\eqn\phiRec{
\phi(P|Q) = {1\over s_P}\!\!\!\sum_{XY=P\atop
AB=Q}\!\!\!\Big(\phi(X|A)\phi(Y|B) - (X\leftrightarrow Y)\Big),\;\;
\phi(i|j) = \d_{ij}\,,\;\phi(\emptyset|B)=\phi(A|\emptyset)\equiv0
}
where $\d_{ij}=1$ if $i=j$ and $0$ otherwise. Note that $\phi(P|Q)=\phi(Q|P)$.
The first instances at multiplicity two are given by
$\phi(12|12)=1/s_{12}=\phi(21|21)$ and $\phi(12|21)=-1/s_{12}$, while at multiplicity three
we have,
$$\eqalignno{
\phi(123|123) &= {1\over s_{123}}\big({1 \over s_{12}}
+ { 1 \over s_{23}}\big),
\quad
\phi(123|132) = - { 1 \over s_{23} s_{123}},
\quad
\phi(123|213) = - { 1 \over s_{12} s_{123}},\cr
\phi(123|321) &=  {1\over s_{123}}\big({1 \over s_{12}}
+ { 1 \over s_{23}}\big),
\quad
\phi(123|231) = - { 1 \over s_{23} s_{123}},
\quad
\phi(123|312) = - { 1 \over s_{12} s_{123}}\,.
}$$
It follows from the antisymmetric deconcatenation in \phiRec\ that
$\phi(P|Q)$ satisfies
the defining symmetry of an {\it alternal mould} \mouldec,
\eqn\alternal{
\phi(P|A\shuffle B) =\phi(A\shuffle B|Q) = 0,\qquad\forall A,B\neq\emptyset\,,
}
where the shuffle product is defined by
\eqn\shuffledef{
\emptyset\shuffle A = A\shuffle\emptyset = A,\qquad
A\shuffle B \equiv a_1(a_2 \ldots a_{n} \shuffle B) + b_1(b_2 \ldots b_{m}
\shuffle A)\,.
}
The identity \alternal\ can be proved by induction \Gauge\ using the linearity
of $\phi(P|Q)$.

For the physics motivation: the construction of the map $\phi(P|Q)$ in \FTlimit\
followed the pioneering work of \DPellis\ where a similar map of planar binary
trees was proposed and used to obtain the tree-level scattering amplitudes of a
theory of bi-adjoint scalars.

%*************************************************
\subsec Planar binary trees and nested bracketings

It is well known that each pb tree can be mapped to a Lie polynomial
in the Free Lie Algebra of the alphabet labelling its
leaves \refs{\garsia,\barcelo}. For example, the trees in \figsingle\ are
mapped to the words $[[1,2],3]=123-213-312+321$ and
$[1,[2,3]]=123-132-231+321$. So let us modify the map \bPT\ of individual
pb trees to also include bracketed words,
\eqn\brPT{
b(P,T) \equiv {1\over s_P}\big[b(P_1,T_1),b(P_2,T_2)\big],\qquad b(i,T)\equiv i\,.
}
The setting is the same as in \bPT:
the tree $T$ is decomposed in terms of its
left $T_1$ and right $T_2$ subtrees with $P_1$ and $P_2$ denoting the
subwords labelling their leaves.
For example, the
two pb trees from \figsingle\ are now mapped to the expressions seen in \figbr.
\ifig\figbr{Example applications of the map $b(P,T)$ defined in \brPT\ for the
two pb trees in $\PBT_3$. The trees are mapped to bracketed words and Mandelstam
variables.}
{\epsfxsize=0.38\hsize\epsfbox{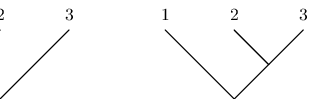}}

The map \brPT\ can be extended to a {\it sum}
over all pb trees
in $\PBT_n$ by deconcatenation: we drop the tree
specification $b(P,T)\rightarrow b(P)$ and evaluate the following recursion:
\eqn\BGbr{
b(P) = {1\over s_P}\sum_{XY=P}[b(X),b(Y)]\,.
}
Just like \BGtree, the deconcatenations generate all trees in $\PBT_n$; but this
time they are dressed with bracketed words and Mandelstam variables. For
example,
\eqn\BGexamples{
b(1) = 1,\quad b(12) = {[1,2]\over s_{12}},
\quad b(123) = {[[1,2],3]\over s_{12}s_{123}}  + {[1,[2,3]]\over s_{23}s_{123}}\,,
}
and see \figBGfour\ for the expression of $b(1234)$.
The same mechanism used in proof of \alternal\ can be used to prove that the
recursion \BGbr\ satisfies,
\eqn\altM{
b(A\shuffle B) = 0\,\quad\forall A,B\neq\emptyset\,.
}
Therefore the additional bracketing structure in the numerators does not
spoil the alternal mould symmetry of the pb trees from \phiRec.

\ifig\figBGfour{The sum over dressed pb trees in $\PBT_4$ generated by the recursion
of $b(1234)$ in \BGbr.}
{\epsfxsize=0.88\hsize\epsfbox{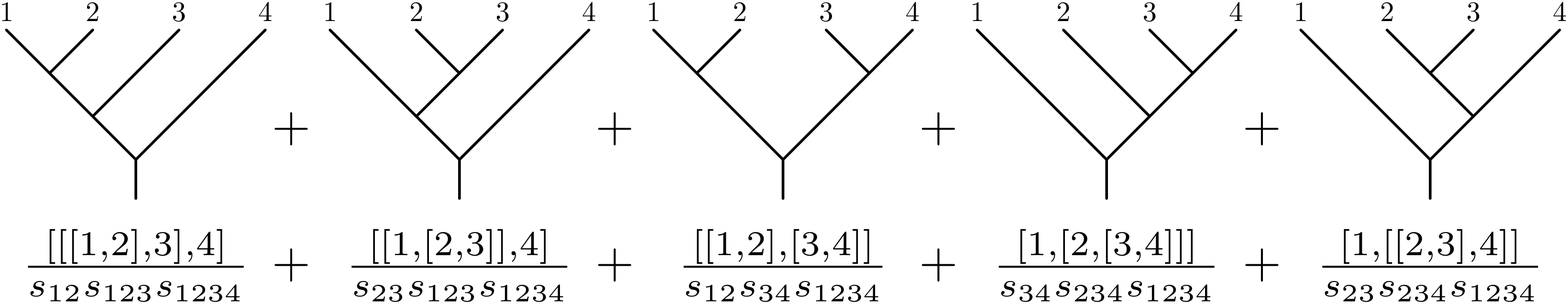}}

%******************************
\newsec Berends--Giele currents

Apart from dressing pb trees with Mandelstam variables and Lie polynomials,
even more structure can be added to them. In doing so, we obtain
objects with direct relevance to the computation of tree-level amplitudes
with the pure spinor formalism. They led to compact expressions
for the amplitudes of both the open superstring
and its field-theory limit \refs{\nptString,\nptMethod}.

%******************************************
\subsec Multiparticle unintegrated vertices

Let us briefly recall the existence of {\it multiparticle unintegrated vertices}
$V_P$; a generalization of the vertex $V_i$ that plays a fundamental role in the
pure spinor formalism \psf. They are defined recursively in a manner described
in \EOMBBs\ but let us focus only on their high-level properties and leave aside
the particularities of their assembly.

The vertices $V_P$ can be characterized by the symmetry relations they satisfy:
the {\it generalized Jacobi identities} as defined in \reutenauer\ (also
referred to as {\it Lie symmetries} in \EOMBBs),
\eqn\genjac{
V_{A\ell(B)C}+V_{B\ell(A)C} = 0\,,\quad A,B\neq\emptyset\,,\quad\forall \, C\,,
}
where $\ell(A)$ is the left-to-right bracketing defined recursively by
\eqn\ellmap{
\ell(123 \ldots n) \equiv \ell(123 \ldots n{-}1)n - n\ell(123 \ldots n{-}1)\,,\quad
\ell(i)\equiv i\,.
}
For example, $V_{1234C} + V_{2143C} + V_{3412C} + V_{4321C} = 0$, for any word
$C$. In addition, the vertices $V_P$ are Grassmann-odd, $V_PV_Q = - V_Q V_P$.
Moreover, with the understanding that a word $P$ inside $V_{ \ldots}$ stands for
$\ell(P)$ and using the notation $V_{\ell(123 \ldots n)} \equiv V_{123 \ldots
n}$ it is always possible to rewrite arbitrary bracketings within $V_{\ldots}$
in terms of $V_P$ with unbracketed $P$, using Baker's identity \reutenauer. For
example, $V_{[A,B]} = V_{A\ell(B)}$ implies that $V_{[1,[[2,3],4]]} = V_{1234} -
V_{1324} - V_{1423} + V_{1432}$. In particular, one can always fix the first
letter of $V_P$ using
\eqn\basisV{
V_{AiB}= -
%(-1)^{\d_{|A|,0}}
V_{i\ell(A)B}\,,\qquad A\neq\emptyset\,.
}
Therefore the number of independent $V_P$ at multiplicity $n$ is
$(n{-}1)!$, in agreement with the dimension of {\it multilinear} Lie
polynomials\foot{A Lie polynomial is called {\it multilinear} when its
words are restricted to be permutations.}
as stated in section 5.6.2 of \reutenauer.

%**********************************************************
\subsec Planar binary trees and nested bracketings of $V_P$

The discussions above can be understood in the context of Free Lie Algebras.
However, in the setting of scattering amplitudes with pure spinor methods, the
bracketed words in the definition \BGbr\ are replaced by vertices $V_{ \ldots}$
with corresponding bracket structure. For example, the word expansions in
\BGexamples\ become {\it Berends--Giele currents} \EOMBBs
\eqn\BGexamplestwo{
M_1 = V_1,\quad M_{12} = {V_{[1,2]}\over s_{12}},\quad
M_{123} = {V_{[[1,2],3]}\over s_{12}s_{123}}  + {V_{[1,[2,3]]}\over
s_{23}s_{123}}\,.
}
Obviously, the alternal
mould property continues to hold and we get
\eqn\altM{
M_{A\shuffle B} = 0\,,\qquad
M_{AiB} = (-1)^{|A|} M_{i(\tilde A \shuffle B)}\,,\quad\forall A,B\neq\emptyset\,,
}
where the second identity (henceforth called Schocker's identity)
was proven in \BGschocker\ and implies that $M_P$ for words of length $n$
admits a $(n{-}1)!$ dimensional basis.

Alternatively,
the Berends--Giele currents
$M_P$ can be defined via a product of $\phi(A|B)$ and {\it unbracketed} superfields $V_B$ \FTlimit:
\eqn\BGdeftwo{
M_{A} = \sum_{B} {1\over|B|}\phi(A|B)V_{B} = \sum_{C} \phi(A|iC)V_{iC}\,.
}
The alternal mould symmetry of $\phi(A|B)$ identifies the sum over $B$ as that
of a Lie polynomial \Ree. This means that the sum over the $|B|!$ permutations
are reduced to a sum over $(|B|{-}1)!$ cyclic permutations of the form $B=iC$,
cancelling the overall factor $1/|B|$. Assuming the equivalence between the
definitions \BGbr\ and \BGdeftwo\ allows to infer a different representation for
$\phi(P|Q)$ as compared to \phiRec,
\eqn\altPhi{
\phi(P|Q) = \langle P, b(Q)\rangle\,,
}
where $\langle A,B\rangle = \d_{A,B}$ is the scalar product of words
and $\d_{A,B}=1$ if $A{=}B$ and $0$ otherwise.
For example, the expansion in  \BGexamples\ together
with $\langle 213, [[1,2],3]\rangle = -1$
and $\langle 213, [1,[2,3]]=0$ implies
\eqn\exbG{
\phi(213|123) = \langle 213, b(123)\rangle = - {1\over s_{12}s_{123}}\,.
}
After defining
\eqn\Phicap{
\Phi(A|B)_i \equiv \phi(iA|iB)
}
the definition \BGdeftwo\ can be rewritten in the form of the so-called {\it BG-map\/}:
\eqn\BGiA{
M_{iA} = \sum_{B} \Phi(A|B)_iV_{iB}\,.
}
The first few examples of \Phicap\ are
$\Phi(2|2)_1 = {1\over s_{12}}$ and
\eqnn\twoPhi
$$\eqalignno{
\Phi(23|23)_1 &= {1 \over s_{12} s_{123}} + { 1 \over s_{23} s_{123}},
\quad
\Phi(23|32)_1 = - { 1 \over s_{23} s_{123}}, &\twoPhi\cr
\Phi(32|32)_1 &= {1 \over s_{13} s_{123}} + { 1 \over s_{23} s_{123}},
\quad
\Phi(32|23)_1 = - { 1 \over s_{23} s_{123}},\cr
}$$
The above definition motivates the following question: Can we obtain a relation
analogous to \BGiA\ where $V_{iA}$ is written in terms of $M_{iB}$? In other words,
can we invert $\Phi(A|B)_i$?

%**************************************************
\subsec The KLT matrix as the inverse of the BG map

We have seen that the alternal mould $\phi(A|B)$ can be used to map the
superfields $V_B$ (with unbracketed words) into the Berends--Giele currents
$M_A$; we are now going to consider its inverse map. We will encounter a
fascinating object called the KLT matrix $S(A|B)_i$ whose origins date back to
the 80s when a relation between amplitudes of closed and open strings was found
\KLTref. The precise relation of \KLTref\ was subsequently formulated in terms
of a KLT matrix in the field-theory limit in \oldMomKer\ and later given its
full string-theory version in \MomKer, with a slight reformulation in \Polylogs.
In the following we will utilize a recent recursive definition given in \NLSM.

More precisely, let $S(A|B)_i$ for words $A$, $B$ and letter $i$ denote
a symmetric matrix that vanishes if $A$ is not a permutation of $B$ and
otherwise given by
\eqn\KLT{
S(P,j|Q,j,R)_i \equiv (k_{iQ}\cdot k_j)S(P|Q,R)_i,\quad
S(\emptyset|\emptyset)_i \equiv 1\,,\quad |Q|+|R|=|P|\,.
}
For an example application of the recursion \KLT\ consider the following sequence:
$S(243|432)_1 = (k_{14}\cdot k_3)S(24|42)_1$,
$S(24|42)_1=(k_1\cdot k_4)S(2|2)_1$ and $S(2|2)_1 = (k_1\cdot k_2)$. Therefore
$S(243|432)_1=(k_{14}\cdot k_3)(k_1\cdot k_4) (k_1\cdot k_2)$.
As an additional example, the entries of the symmetric $2!\times2!$ matrix composed
of permutations from $A,B{=}23$ with $i{=}1$ are:
\eqn\twoS{
S(23|23)_1 = s_{12}(s_{13}+s_{23})\,,\quad
S(23|32)_1 = s_{12}s_{13}\,,\quad
S(32|32)_1 = s_{13}(s_{12} + s_{23})\,.
}
As one can easily check from the examples above, we have $S(2|2)_1\Phi(2|2)_1= 1$ as well as
\eqnn\Stwoex
$$\eqalignno{
S(23|23)_1\Phi(23|23)_1
+ S(23|32)_1\Phi(32|23)_1 &= 1 &\Stwoex\cr
S(23|23)_1\Phi(23|32)_1
+ S(23|32)_1\Phi(32|32)_1 &= 0\,.
}$$
Following the arguments of \DPellis\ in the context of scattering amplitudes and
its reformulation in terms of $\phi(A|B)$ from \FTlimit,
there is a strong expectation that this must be true in general:
\eqn\InvSPhi{
\sum_C S(A|C)_i \Phi(C|B)_i = \d_{A,B}\,.
}
The sum over $\sum_{C}$ instructs to sum
over all words $C$ but the condition that $S(A|C)_i$ is zero if
$C$ is not a permutation of $A$ leads to a finite sum. In the subsequent
discussions we will obtain a more general relation that reduces to
\InvSPhi\ in a particular case, leading to an alternative avenue
to prove it.

%***************************************
\subsubsec KLT matrix and labelled trees

If $s_{ij}$ is depicted as the edge between the vertices labelled $i$ and $j$,
the entries of the KLT matrix $S(A|B)_i$ for words $A,B$ of length $n$
generate all rooted labelled trees with $n+1$ vertices and $n$ edges, whose
total number is $(n+1)^{n-1}$ \stanley.

To see this, note that the recursion \KLT\ removes the letter $j$ to generate an
edge $k_R\cdot k_j$ at each iteration for a total of $n$ edges; in particular,
$(k_i\cdot k_j)^2$ is never generated. In addition, we see from
\KLT\ that there is always a path from any given vertex to the vertex~$i$.
Therefore the recursion \KLT\ gives rise to rooted labelled trees. Conversely,
if a given tree with $n$ vertices appears in $S(A|B)_i$ then a new tree obtained
by appending an edge $s_{j(n+1)}$ where $j=1, \ldots n$ is necessarily contained
in $S(A(n{+}1)|B(n{+}1))_i=k_{iB}\cdot k_{n+1}S(A|B)_i$. Relabelling the vertices if
necessary, the new vertex $n{+}1$ can always be chosen to appear in a single edge
of the new tree where the argument above applies, finishing the proof.

For example, the tree representation of $S(234|423)_1=s_{12}s_{13}s_{14} +
s_{14}s_{12}s_{23}$ is:
\medskip
\centerline{{\epsfxsize=0.70\hsize\epsfbox{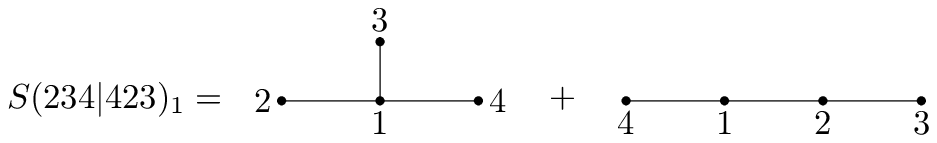}}}
\medskip
\noindent For another application, note that a symmetric $3!\times 3!$ matrix
naively contains $21$ elements. By the above proof there are
only $4^2=16$ cubic monomials of $s_{jk}$ in the permutations of
$S(234|234)_1$, so there must be five relations among them.
An explicit search yields:
\eqnn\relat
$$\eqalignno{
S(432|234)_1 &= S(342|243)_1,\quad S(423|324)_1 = S(342|243)_1\,, &\relat\cr
S(432|342)_1 &= S(423|342)_1 - S(342|243)_1 + S(432|324)_1\,,\cr
S(432|243)_1 &= S(342|243)_1 - S(423|234)_1 + S(423|243)_1\,,\cr
S(342|243)_1 &= S(324|243)_1 - S(324|234)_1 + S(342|234)_1\,.
}$$
By the same token, any other choice of the fixed letter $i$ can
be expanded in terms of a basis where $i=1$.

It is amusing to note that $(n{+}1)^{n{-}1}$ is also the number of ways to
factorize the cycle $(12 \ldots (n{+}1))$ as a product of transpositions
\cycfac. For example, the factorizations of the cycle $(123)$ corresponding to
the three labelled trees with three vertices are $(123) {=} (12)(23) {=}
(13)(12) {=} (23)(13)$ (see e.g. exercise 5.47 of \stanley).

%*****************************
\subsec An extended KLT matrix

Given the conjectural relation \InvSPhi\ and the definition \BGdeftwo\ it is
not difficult to see that the KLT matrix inverts \BGdeftwo\
leading to \Polylogs,
\eqn\MToV{
V_{iA} = \sum_{B}S(A|B)_i M_{iB}\,.
}
Unlike the definition \BGdeftwo\ where $\phi(\cdot|\cdot)$ can be used for
arbitrary permutations of both words (thereby manifesting the alternal mould
symmetry of $M_A$) the relation \MToV\ fixes the first letter in the left-hand
side to be $i$. This is unsatisfactory since we know from \genjac\ that $V_P$
satisfies generalized Jacobi identities when all permutations of $P$ are
considered but \MToV\ manifestly defines only a cyclic orbit of $P\equiv iA$.
Although nothing prevents using different choices of $i$ on demand, a more
general definition akin to \BGdeftwo\ is desirable.

The lack of manifest Lie symmetry in \MToV\ is due to the
definition of $S(A|B)_i$, which obviously singles out $i$. Therefore
the quest is to find a general definition for the KLT matrix
without this restriction. Luckily, this can be done with the following
arguments.

One can explicitly verify in examples
that when $S(P|Q)_i$ is considered as a function of permutations
of the extended words $iP,iQ$ it satisfies the following constraint:
\eqn\conjS{
S(AjB|RjS)_i = S(\ell(iA)B\,|\,\ell(iR)S)_j\,.
}
For example, $S(314|314)_2 = S(234|234)_1 - 2S(234|324)_1 + S(324|324)_1$.

The identities \basisV\ and \conjS\ suggest the
existence of a more general KLT matrix $S(A|B)$ with no fixed letter $i$ that
reduces to \KLT\ when the first letters of $A,B$ coincide
\eqn\SiP{
S(iP|iQ) \equiv S(P|Q)_i\,.
}
To see this, note that if this extended
matrix $S(A|B)$ satisfies Lie symmetries
in both row and columns, namely (note $|A|=|B|$),
\eqn\genLieKLT{
S(\ell(A)|B)=S(A|\ell(B))=|A|S(A|B)\,,
}
then the identity \conjS\ is explained,  according to \basisV,
as a ``change of basis'' of multilinear Lie polynomials whose first letter
is fixed to be $i$ or $j$.
This observation can be exploited to define
an {\it extended KLT matrix\/} $S(A|B)$ as a multilinear Lie polynomial
in both (arbitrary) words $A$ and $B$ as follows:
Without loss of generality we write $A=PiQ$ and $B=RiS$
and move the letter $i$ to the front via \basisV\ to obtain,
\eqn\genKLT{
S(P,i,Q|R,i,S)
\equiv (-1)^{\d_{|P|,0}}(-1)^{\d_{|R|,0}}S(\ell(P)Q|\ell(R)S)_i\,,
}
where the sign factors account for the possibility of either $P$ or $Q$ being
the empty word of length zero\foot{The sign in $PiQ=-i\ell(P)Q$
is negative only if $P\neq\emptyset$ as there is nothing to do in case $P=\emptyset$.
We set $\ell(\emptyset)\equiv\emptyset$ and
write $PiQ=-(-1)^{\d_{|P|,0}}i\ell(P)Q$ to extend its validity also when $P=\emptyset$.
}.
One can now check that \genKLT\ and \SiP\ imply the identity \conjS:
\eqn\conjSequiv{
S(AjB|RjS)_i \equiv S(iAjB|iRjS) \equiv S(j\ell(iA)B|j\ell(iR)S)
\equiv S(\ell(iA)B|\ell(iR)S)_j\,.
}
For example, the entries of the extended KLT matrix for
sample permutations of
$A,B{=}123$ are given by
\eqnn\sampleP
$$\eqalignno{
S(213|213) &= S(123|123) = S(23|23)_1\,, &\sampleP\cr
S(213|231) &= S(123|123) - S(123|132) = S(23|23)_1 - S(23|32)_1\,,\cr
}$$
which are reduced to the cases that can be computed by \KLT.
Note that {\it any} common letter $i$ between $A$ and $B$ can be chosen
in \genKLT\ because \conjS\ guarantees their equality. For example,
$S(213|321)$ with $i{=}1$ leads to
$S(213|321) = S(23|32)_1 - S(23|23)_1 = - s_{12}s_{23}$
whereas $i{=}2$ gives $S(213|321) = -S(13|31)_2 = - s_{12}s_{23}$.

%*********************************************************
\subsubsec Extended KLT matrix as an inverse to the BG map

Using the general matrix $S(A|B)$ given in \genKLT\ the definition
\MToV\ can be promoted to
\eqn\genMToV{
V_A = \sum_{B} {1\over |B|} S(A|B)M_B\,.
}
Given \genLieKLT, the definition \genMToV\ manifests
the Lie symmetries of $V_A$ and accomplishes the task
initiated in the last subsection.
But now the compatibility between \MToV\ and \BGdeftwo\ suggests
a relation between $S(A|B)$ and $\phi(A|B)$. Experimentally one finds,
\eqn\genInv{
\sum_B S(A|B)\phi(B|C) = \langle \ell(A), \rho(C)\rangle\,,
}
where $\ell(A)$ is defined in \ellmap\ while $\rho(C)$ is given by \reutenauer
\eqn\rhomap{
\rho(123 \ldots n) \equiv 1\rho(23 \ldots n) - n\rho(123 \ldots n{-}1),
\qquad\rho(i)\equiv i\,,
}
and $\langle P,Q\rangle = \d_{P,Q}$ denotes the scalar product of words.
From the fact that $\rho$ and $\tilde\ell$ are adjoint maps w.r.t the
scalar product together
with $\tilde\ell\ell(A)=|A|\ell(A)$ and $\tilde\rho\rho(C)=|C|\rho(C)$
\reutenauer,
it follows that the RHS of \genInv\ can be written as
$|C|\langle A, \rho(C)\rangle = |A|\langle \ell(A), C\rangle$
and that it is not positive definite.

We can now demonstrate that \genInv\ implies
the conjectural relation \InvSPhi.

\proclaim Proposition 1. If
$\sum_C \phi(A|C)S(C|B) = \langle \rho(A), \ell(B)\rangle$
then \InvSPhi\ is true,
\eqn\true{
\sum_R \Phi(P|R)_i S(R|Q)_i = \d_{P,Q}\,.
}\par
\noindent{\it Proof:\/} The shuffle symmetry
of $\phi(A|C)$ makes the sum over all $C$ reduce to a cyclic
subset of permutations with an overall
factor $|C|{=}|B|$. Setting $C{\equiv} iC'$, we get
\eqn\tmpI{
\langle \rho(A), \ell(B)\rangle =
\sum_C \phi(A|C)S(C|B) = \sum_{C'} |B| \phi(A|iC')S(iC'|B)\,.
}
Since $\rho$ is the adjoint of the right-to-left
bracketing map $r$ and $r(P){=}|P|P$ if $P$ is a Lie
polynomial (which $\ell(B)$ certainly is) \reutenauer;
$\langle \rho(A),\ell(B)\rangle {=} |B|\langle A,\ell(B)\rangle$.
Therefore we get
$\sum_{C'} \phi(A|iC')S(iC'|B) = \langle A,\ell(B)\rangle$.
Choosing $A{\equiv} iA'$
and $B{\equiv} iB'$ leads to
\eqn\finalq{
\sum_{C'}\Phi(A'|C')_i S(C'|B')_i = \langle iA', \ell(iB')\rangle =
\d_{A',B'}\,,
}
where we used \Phicap\ and \SiP. To prove the last equality,
note that the identity
$\ell(aP) = \sum_{X\shuffle Y=P}(-1)^{|X|} \tilde XaY$ \patras\ implies
$\ell(iB') = iB' + \sum RiS$ with $R\neq\emptyset$ and
therefore
$\langle iA',\ell(iB')\rangle = \d_{A',B'}$. Renaming $A'{=}P, B'{=}Q$ and
$C'{=}R$ finishes the proof \qed.

In addition,
the extended KLT matrix \genKLT\ satisfies
\eqn\eKLTob{
S(A|b(C)) = \langle\ell(A),\rho(C)\rangle\,,
}
where the $b(C)$ map is defined in \BGbr.
This follows from plugging in $\phi(B|C) = \langle B, b(C)\rangle$ into the relation
\genInv\ and using the definition of the scalar product $\langle
P,Q\rangle=\d_{P,Q}$.

For example, on the one hand
$\langle \ell(123),\rho(123)\rangle = 3$ while on the other
hand the expansion \BGexamples\ of $b(123)$ leads to
\eqnn\SAbC
$$\eqalignno{
S(123|b(123)) &= {1\over
s_{12}s_{123}}\Big(S(123|123)-S(123|213)-S(123|312)+S(123|321) \Big)\cr
&+ {1\over s_{23}s_{123}}\Big(S(123|123)-S(123|132)-S(123|231)+S(123|321)
\Big)\cr
&={1\over s_{12}s_{123}}\Big(3S(23|23)_1\Big)
+ {1\over s_{23}s_{123}}\Big(3S(23|23)_1-3S(23|32)_1\Big)\cr
&= 3\,,&\SAbC
}$$
where in the last step we used the expressions \twoS.

We have defined the extended KLT matrix by reducing its
permutations to the standard KLT matrix via \SiP\
where the recursive algorithm \KLT\ can be applied.
After introducing the so-called S-map defined in \EOMBBs,
we will obtain a direct formula to compute the entries of the
extended KLT matrix without recoursing to its old definition \KLT.

%****************************************************************
\subsec The S-map and bracketed numerators of planar binary trees

It is easy to see from \BGbr\ that $[1,2] = s_{12}b(12)$, but already
at the next order
\eqnn\notr
$$\eqalignno{
[[1,2],3] &= s_{12}s_{23} b(123) - s_{12}s_{13}b(213)\,,&\notr\cr
[1,[2,3]] &= s_{12}s_{23} b(123) - s_{13}s_{23}b(132)\,,
}$$
an interplay between the Jacobi identity %$[1,[2,3]]= [[1,2],3] - [[1,3],2]$ 
and $s_{123}=s_{12} + s_{13} +s_{23}$ is needed
to cancel all denominators.
The identities in \notr\ motivate the search for a general
procedure that generates the expansions on the right-hand side for any
given numerator
of an arbitrary pb tree. Surprisingly, an algorithm
discovered in \EOMBBs\ (while in pursuit of other objectives) can be used
to achieve precisely that.

The algorithm is based on the so-called S-map $\{A,B\}$ between words $A$ and $B$
defined as follows:
\eqn\Smap{
\{A,B\} \equiv (-1)^{|B|+1} \rho(A)\otimes^s \tilde\rho(B)\,,
}
where $\rho(B)$ is defined in \rhomap\ and $\tilde\rho(B)$ denotes
its reversal while
$\otimes^s$ denotes a {\it weighted} concatenation product between words
\eqn\wconcdef{
Ai \otimes^s jB \equiv s_{ij} AijB\,.
}
For example, given $\rho(12)= 12 - 21$ and
$\tilde\rho(3) = 3$, we get $\{12,3\}=s_{23}123 - s_{13}213$.

The claim based on experimental data is that the map $b(A)$
defined in \BGbr\ acting on a nested
application of the S-map gives rise to a pb tree numerator with corresponding
bracketed structure
\eqn\nestS{
[[ \ldots [i,j],k] \ldots] = b\big(\{\{ \ldots \{i,j\},k\} \ldots\}\big)\,.
}
For example, the relations in \notr\ are obtained from
\eqn\notrEx{
[[1,2],3] = b\big(\{\{1,2\},3\}\big)\,,\qquad
[1,[2,3]] = b\big(\{1,\{2,3\}\}\big)\,,
}
while $[[1,2],[3,4]] = b\big(\{\{1,2\},\{3,4\}\}\big)$ implies
\eqn\notrt{
[[1,2],[3,4]]
%       + M(1,2,3,4)*s(1,2)*s(2,3)*s(3,4)
%       - M(1,2,4,3)*s(1,2)*s(2,4)*s(3,4)
%       - M(2,1,3,4)*s(1,2)*s(1,3)*s(3,4)
%       + M(2,1,4,3)*s(1,2)*s(1,4)*s(3,4)
=s_{12}s_{34}\Big(s_{23}b(1234)
- s_{24}b(1243)
- s_{13}b(2134)
+ s_{14}b(2143)\Big)\,,
}
which can be verified with some effort.
The examples above suggest that it is convenient to define
a nested left-to-right S-map by
\eqn\sigmamap{
\sigma(123 \ldots n) \equiv \big\{ \sigma(12 \ldots n{-}1), n\big\},\qquad
\sigma(i) \equiv i\,,
}
for example $\s(123) = \{\{1,2\},3\}=s_{12}(s_{23}\,123 - s_{13}213)$.
Using \sigmamap\ the conjectural relation \nestS\ implies
\eqn\ellsigma{
\ell(A) = b(\sigma(A))\,,
}
which can be verified at high multiplicities.
The S-map was originally defined in \EOMBBs\ in terms of
Berends--Giele currents $M_A$ using the notation $M_{S[A,B]} \equiv
M_{\{A,B\}}$ as follows
\eqn\SmapOrig{
M_{S[A,B]} \equiv \sum_{i=1}^{|A|} \sum_{j=1}^{|B|} (-1)^{i-j+|A|-1}
s_{a_i b_j} M_{(a_1 a_2\ldots a_{i-1} \shuffle
a_{|A|} a_{|A|-1}\ldots a_{i+1})a_ib_j
(b_{j-1}\ldots b_2 b_1 \shuffle b_{j+1} \ldots b_{|B|})}
}
for $A=a_1 a_2\ldots a_{|A|}$ and $B=b_1 b_2\ldots b_{|B|}$. The equivalence
between the definition \Smap\ and \SmapOrig\ follows from the following
general identity of words
\eqn\rhoasshuffle{
\rho(A) = \sum_{XjY=A}(-1)^{|Y|}(X\shuffle\tilde Y)j\,,
}
which can be proven by induction using the definitions \rhomap\ and \shuffledef.

%*****************************************************
\subsubsec Grafting of trees and an algebra for $b(A)$

Recall that the {\it grafting} operation on pb trees $T_1$ and $T_2$ is denoted
$T_1 \vee T_2$ and gives rise to a new pb tree where the roots of $T_1$ and
$T_2$ are glued together \refs{\vallette,\loday}. For example,
\bigskip
\centerline{{\epsfxsize=0.50\hsize\epsfbox{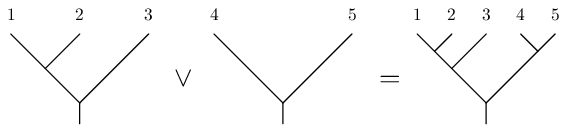}}}
\medskip
\noindent Given that each pb tree corresponds to a Lie polynomial composed of bracketed
words, the grafting of trees induces a map between Lie polynomials. For
example, the above grafting implies $[[1,2],3] \vee [4,5] = [[[1,2],3],[4,5]]$.
An interesting observation (confirmed by explicit examples)
is that the S-map seems to capture the effect of
grafting all the trees and associated Mandelstam variables in $b(A)\vee
b(B)$. More precisely,
\eqn\bAlgebra{
b(A)\vee b(B) = b\big(\{A,B\}\big)\,,
}
giving rise to an algebraic structure of $b(A)$ trees.

\ifig\figgraft{The graphical depiction of the grafting $b(123)\vee
b(45)$. The right-hand side
can be expressed in terms of the S-map as $b\big(\{123,45\}\big)$.}
{\epsfxsize=0.94\hsize\epsfbox{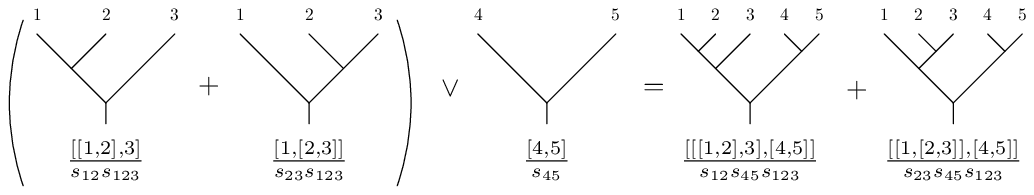}}

A non-trivial example of
the above algebra is obtained by considering the grafting of $b(123)\vee b(45)$.
On the one hand we get (see \figgraft),
\eqn\lhsalg{
b(123)\vee b(45) = {[[[1,2],3],[4,5]]\over s_{12}s_{45}s_{123}}
+ {[[1,[2,3]],[4,5]]\over s_{23}s_{45}s_{123}}\,.
}
On the other hand, the S-map between the words $A{=}123$ and $B{=}45$ is given by
\eqnn\rhsalgtmp
$$\eqalignno{
% smap of {123,45} =
       % + Word(1,2,3,4,5)*s(3,4)
       % - Word(1,2,3,5,4)*s(3,5)
       % - Word(1,3,2,4,5)*s(2,4)
       % + Word(1,3,2,5,4)*s(2,5)
       % - Word(3,1,2,4,5)*s(2,4)
       % + Word(3,1,2,5,4)*s(2,5)
       % + Word(3,2,1,4,5)*s(1,4)
       % - Word(3,2,1,5,4)*s(1,5)
\{123,45\} &=
         s_{34}12345
        - s_{35}12354
        - s_{24}13245
        + s_{25}13254
        - s_{24}31245 &\rhsalgtmp\cr
        &\quad{}+ s_{25}31254
        + s_{14}32145
        - s_{15}32154\,.
}$$
A long calculation using the $b(A)$ map \BGbr\ together with the
linearity condition \funlin\ yields
\eqn\rhsalg{
b\big(\{123,45\}\big) = {[[[1,2],3],[4,5]]\over s_{12}s_{45}s_{123}}
+ {[[1,[2,3]],[4,5]]\over s_{23}s_{45}s_{123}}\,.
}
Therefore we see from \lhsalg\ and \rhsalg\ that $b(123)\vee b(45) = b\big(\{123,45\}\big)$.

%************************************************************
\subsec An alternative definition for the extended KLT matrix

The introduction of the S-map in the last subsection suggests that a direct
definition for the extended KLT matrix which does not rely on the standard KLT
matrix as in \genKLT\ may be possible. The reasoning is similar to the one that
led us to the relation \altPhi\ and exploits the fact that $V_P$ can be obtained
from two different ways. Consider the example of $V_{123}$. On the one hand it
can be obtained from the definition \genMToV,
\eqn\tmpIa{
V_{123} = \sum_B {1\over 3}S(123|B)M_B = S(123|123)M_{123} + S(123|132)M_{132}\,.
}
On the other hand it follows from the S-map on the Berends--Giele currents
\eqn\tmpIb{
V_{123} = s_{12}s_{23}M_{123} - s_{12}s_{13}M_{213}
= (s_{12}s_{23} + s_{12}s_{13})M_{123} + s_{12}s_{13}M_{132}\,,
}
by exploiting the relation of $[[1,2],3]\rightarrow V_{123}$ and
$b(123)\rightarrow M_{123}$. In the last equality we used Schocker's identity \altM\
to rewrite $M_B$ in a basis of $M_{1B'}$. Comparing \tmpIa\ and \tmpIb\ leads
to the expressions for $S(123|123)$ and $S(123|132)$. Note that both definitions
do not necessarily rely on any particular basis of $M_{iA'}$ and are valid in
general. Thus the S-map gives rise to a recipe for unlocking
the expressions of $S(A|B)$ for arbitrary permutations. In fact, the
considerations to be given below lead to the following proposal,
\eqn\altS{
S(P|Q) \equiv \langle \ell(P), \s(Q)\rangle.
}
For example, from $\ell(123) = 123 - 213 - 312 + 321$ and
$\s(123) = s_{12}s_{23}123 - s_{12}s_{13}213$ one immediately gets
$S(123|123) = s_{12}(s_{13} + s_{23})$. Higher-multiplicity examples
are similarly verified.

After the experimental observation that $\sigma(AjD) = - k_A\cdot k_j\,
j\s(AD) + \ldots$ where the omitted terms do not
contain the letter $j$ at the last position when\foot{
The definition of \sigmamap\ implies that the last letter in $\s(A)$
is always the last letter of $A$.} $D\neq\emptyset$ one can see that \altS\ reduces to
the recursion \KLT\ for $S(Aj|CjD)_i$,
\eqnn\reduce
$$\eqalignno{
S(Aj|CjD)_i &= \langle \ell(iAj), \s(iCjD)\rangle =
- k_{iC}\cdot k_j \langle \ell(iA)j - j\ell(iA),  j\sigma(iCD) + \cdots \rangle\cr
& = k_{iC}\cdot k_j \langle \ell(iA),  \sigma(iCD) \rangle
 = k_{iC}\cdot k_j S(A|CD)_i\,,\qquad D\neq\emptyset\,. &\reduce
}$$
In addition, it is straightforward to see using \altS\ and \altPhi\ the
proof of \genInv\ reduces to showing that
$\langle \ell(A), \sigma(b(C))-\rho(C)\rangle=0$, which has been verified
to high multiplicity.

A few words about the proposal \altS\ are in order. In the left-hand side of
\tmpIb, one can interpret $V_A$ as an ``ordinary'' word $A$ but $M_A$ in the
right-hand side does not admit such an interpretation as they satisfy the
shuffle symmetries \altM. This means that one cannot read off the coefficients
$S(A|B)$ from the standard scalar product of words in a similar fashion as in
\altPhi; a different prescription is needed. The simplest trial is as follows:
If $S_P$ satisfies shuffle symmetries as $S_{A\shuffle B}=0$ then we define
$\langle A, S_P \rangle_s$ between an ordinary word $A\equiv iB$ and $S_P\equiv
S_{CiD}$ as
\eqn\schprod{
\langle A, S_P\rangle_s = \langle iB, S_{CiD}\rangle_s
\equiv \langle iB, (-1)^{|C|}i(\tilde C\shuffle D)\rangle
= (-1)^{|C|}\langle B, \tilde C\shuffle D\rangle\,,
}
in terms of the standard scalar product.
Using the above definition and the identity
$\ell(iB) = \sum_{X\shuffle Y=B}(-1)^{|X|} \tilde XiY$ \patras\
one can show that
\eqn\shp{
\langle A, S_P\rangle_s = \langle \ell(A), P\rangle\,.
}
Now define the shuffle-extension $\s_s(A)$ of $\s(A)$ by mapping the words in
the right-hand side of \sigmamap\ to shuffle-satisfying objects according to
$P\rightarrow S_P$. For example, $\sigma(12) = s_{12}12$ while $\sigma_s(12) =
s_{12}S_{12}$. Similarly $\sigma(123) = s_{12}s_{23}123 - s_{13}s_{23}213$ while
$\sigma_s(123) = s_{12}s_{23}S_{123} - s_{13}s_{23}S_{213}$. Using the
definitions above and the intuition gained from the example \tmpIb\ suggests
that $S(P|Q)$ can be extracted using the shuffle-aware scalar product \schprod\
as $S(P|Q) \equiv \langle P, \sigma_s(Q)\rangle_s$, which is {\it
experimentally} checked to be correct. The identity \shp\ then leads to
\altS\ in terms of the standard scalar product.

%*************************************************
\subsec Deconcatenation of Berends--Giele currents

In the pure spinor formalism there is a nilpotent Grassmann-odd BRST operator
$Q^2=0$ \psf. It was argued in \EOMBBs\ that a very interesting pattern arises
in the computation of $QV_P$:
\eqnn\exampOne
$$\eqalignno{
QV_1 &= 0\,,\qquad QV_{12} = (k_1\cdot k_2)V_1V_2 &\exampOne\cr
QV_{123} & = (k_1\cdot k_2)\big[V_1 V_{23} + V_{13}V_2\big]
+ (k_{12}\cdot k_3) V_{12} V_3 \cr
Q V_{1234} &=(k_1\cdot k_2)\bigl[V_1V_{234}
+ V_{13} V_{24} + V_{14} V_{23} +  V_{134} V_2 \bigr]\cr
&\quad{} + (k_{12}\cdot k_3)\bigl[V_{12} V_{34}
+  V_{124} V_3\bigr] + (k_{123}\cdot k_4) V_{123} V_4\,.
}$$
A non-trivial consistency check of the above identities consists in checking
that the right-hand side preserves the symmetries \genjac\ of
the left-hand side (we use that $V_P$ is
Grassmann-odd). It turns out that these identities admit the following
generalization, as suggested from the equations of motion of the superfields
defined recursively in \EOMBBs\ ($k_\emptyset\equiv0$)
\eqn\QVP{
QV_{P} = \sum_{P=XjY\atop Y=R\shuffle S}(k_{X} \cdot k_j)\,
V_{XR}V_{jS}\,,
}
where $Y=R\shuffle S$ represents the {\it deshuffle} of $Y$ into the words $R$
and $S$. An algorithmic procedure to obtain the pairs $(R,S)$ in the sum in
\QVP\ follows from the map $\d_2(Y)$ \reutenauer,
\eqn\deltamap{
\d_2(Y) = \sum_{R,S}\langle Y, R\shuffle S\rangle\, R\otimes S\,.
}
For example $\d_2(123)=
\emptyset\otimes123 + 1\otimes23 + 2\otimes13 + 12\otimes3 + 3\otimes12
+ 13\otimes2
+ 23\otimes1
+ 123\otimes\emptyset$.
It turns out that the identity \QVP\ gives rise to a beautiful
deconcatenation formula for $QM_P$ using the definition
\BGdeftwo. For example, note that the $1/s_{12}$ factor in the
definition of $M_{12}$ cancels the numerator
from $QV_{12}=s_{12}V_1V_2$ so
that $QM_{12} = M_1 M_2$. Similarly and rather surprisingly one
finds the precise cancellations of numerators and denominators to
get $QM_{123} = M_1 M_{23} + M_{12} M_3$. It has been experimentally
checked to high multiplicities that the following fascinating
identity holds true ($M_\emptyset\equiv0$)
\eqn\QM{
QM_P = \sum_{P=XY} M_X M_Y\,.
}
It would be desirable to {\it prove} \QM\ using the
combinatorial definitions \BGdeftwo\ and \QVP.

The above Berends--Giele currents have been constructed in pursuit of a general
formula for the $n$-point scattering amplitude of super-Yang--Mills at
tree-level using BRST cohomology methods \nptMethod. The formula reads
\eqn\SYM{
A(123 \ldots n) = E_{123 \ldots n{-}1} M_n\,,\qquad E_P\equiv
\sum_{P=XY}M_XM_Y\,,
}
and has been shown in \BGPS\ to be the supersymmetric generalization of the
standard Berends--Giele recursion given in \BerendsME.

\bigskip \noindent{\bf Acknowledgements:}
I want to thank Oliver Schlotterer for collaboration on
related topics and the organizers of CARMA2017 for the kind invitation.
I also acknowledge support by a University Research
Fellowship from the Royal Society.

\listrefs
\bye